\documentclass[11pt,namelimits,sumlimits,a4paper]{amsart}
\usepackage{cite}
\usepackage{comment}

\usepackage{amssymb,amsmath}
\usepackage[mathscr]{eucal}
\usepackage{slashed}

\usepackage[latin1]{inputenc}
\usepackage[T1]{fontenc}
\usepackage[UKenglish]{babel}
\usepackage{amsfonts}
\usepackage{fancyhdr}
\usepackage{amsmath}
\usepackage{amsthm}
\usepackage{amsmath,amscd}
\usepackage{latexsym}
\usepackage{cite}
\usepackage{amssymb,amsmath}

\usepackage[mathscr]{eucal}
\usepackage{slashed}
\usepackage{times,psfrag}
\usepackage{mathtools}
\usepackage{hyperref}
\usepackage{multirow}
\usepackage{longtable}
\usepackage{bm}
\usepackage{fancyhdr}
\usepackage{float}

\usepackage{esint}

\usepackage{layout}

\usepackage{lmodern}
\usepackage[rgb]{xcolor}
\usepackage[draft,author={Lorenzo Foscolo}]{pdfcomment}

\usepackage{enumitem}

\numberwithin{equation}{section}
\newtheorem{theorem}[equation]{Theorem}
\newtheorem*{theorem*}{Theorem}
\newtheorem{lemma}[equation]{Lemma}
\newtheorem{prop}[equation]{Proposition}

\theoremstyle{definition}
\newtheorem{definition}[equation]{Definition}

\theoremstyle{remark}
\newtheorem{remark}[equation]{Remark}
\newtheorem*{remark*}{Remark}

\newcommand{\ie}{\emph{i.e.} }
\newcommand{\eg}{\emph{e.g.} }
\newcommand{\cf}{\emph{cf.} }

\newcommand{\beq}{\begin{equation}}
\newcommand{\eeq}{\end{equation}}
\newcommand{\bea}{\begin{eqnarray}}
\newcommand{\eea}{\end{eqnarray}}

\newcommand{\C}{\mathbb{C}}

\newcommand{\R}{\mathbb{R}}

\newcommand{\CP}{\mathbb{CP}}
\newcommand{\PP}{\mathbb{P}}

\newcommand{\ra}{\rightarrow}

\newcommand{\vol}{\operatorname{Vol}}
\newcommand{\dvol}{\operatorname{dv}}

\newcommand{\Real}{\operatorname{Re}}
\newcommand{\Imag}{\operatorname{Im}}

\newcommand{\Lie}[1]{\mathfrak{#1}}
\newcommand{\F}{\mathbb{F}}

\newcommand{\gtwo}{\textup{$G_2$ }}

\newcommand{\sunitary}[1]{\textup{$SU(#1)$}}

\newcommand{\spin}[1]{\textup{$Spin(#1)$}}

\def\co{\colon\thinspace}

\setlist{leftmargin=*}

\begin{document}

\title{Deformation theory of nearly K\"ahler manifolds}

\author[L.~Foscolo]{Lorenzo~Foscolo}
\address{Mathematics Department, State University of New York at Stony Brook}
\email{lorenzo.foscolo@stonybrook.edu}

\maketitle

\begin{abstract}
Nearly K\"ahler manifolds are the Riemannian $6$--manifolds admitting real Killing spinors. Equivalently, the Riemannian cone over a nearly K\"ahler manifold has holonomy contained in $G_2$. In this paper we study the deformation theory of nearly K\"ahler manifolds, showing that it is obstructed in general. More precisely, we show that the infinitesimal deformations of the homogeneous nearly K\"ahler structure on the flag manifold are all obstructed to second order.
\end{abstract}

\section{Introduction}

A \emph{Killing spinor} on a Riemannian spin manifold $(M^n,g)$ is a spinor $\psi$ such that
\[
\nabla_X\psi = \alpha\, X\cdot\psi,
\]
for some $\alpha\in\C$ and all vector field $X$. Here $ \cdot $ denotes Clifford multiplication. Killing spinors appeared in work of Friedrich \cite{Friedrich} on the first eigenvalue of the Dirac operator. It was shown in \cite{Friedrich} that every manifold with a Killing spinor is Einstein: $\text{Ric}(g)=4(n-1)\alpha^2 g$. In particular, one of three cases must hold: (i) $\alpha\neq 0$ is purely imaginary and $M$ is non-compact; (ii) if $\alpha=0$ then $\psi$ is a parallel spinor and therefore $(M,g)$ has holonomy contained in $SU(\frac{n}{2})$, $Sp(\tfrac{n}{4})$, \gtwo or \spin{7}; (iii) $\alpha \neq 0$ is real: $\psi$ is called a \emph{real Killing spinor} and $M$ (if complete) is compact with finite fundamental group. In the real case one can always assume that $\alpha=\pm\tfrac{1}{2}$ by scaling the metric.

There is in fact a relation between parallel spinors and real Killing spinors: by work of B\"ar \cite{Bar}, the cone over a manifold with a real Killing spinor has a parallel spinor and conversely simply connected manifolds with a real Killing spinor are the cross-sections of Riemannian cones with holonomy contained in $SU(\frac{n}{2})$, $Sp(\tfrac{n}{4})$, \gtwo or \spin{7}, depending on the dimension $n$ and the number of linearly independent real Killing spinors.

\emph{Nearly K\"ahler} manifolds are the $6$--dimensional Riemannian manifolds admitting real Killing spinors. The cone over a nearly K\"ahler manifold has holonomy contained in $G_2$.

\begin{remark}
The name nearly K\"ahler was introduced by Gray \cite{Gray} to denote a special class of almost Hermitian manifolds in every even dimension. What we call here nearly K\"ahler manifolds are often referred to as \emph{strict} nearly K\"ahler manifolds of dimension $6$. The terminology \emph{Gray manifolds} has also been used, \cf \cite[Definition 4.1]{Moroianu:Nagy:Semmelmann}.
\end{remark}

Despite the spinorial point of view will play a role in this paper, we prefer to relate the holonomy reduction of the cone $C(M)$ over a nearly K\"ahler manifold to the existence of a closed and co-closed \emph{stable} $3$--form rather than to the existence of a parallel spinor. From this point of view a nearly K\"ahler structure is an \sunitary{3} structure with special torsion: a pair of differential forms $(\omega,\Omega)$, where $\omega$ is a non-degenerate $2$--form and $\Omega$ is a complex volume form satisfying appropriate algebraic compatibility conditions and the first order PDE system
\[
d\omega=3\Real\Omega, \qquad d\Imag\Omega = -2\omega^2.
\]

There are currently only six known examples of simply connected nearly K\"ahler manifolds. Four of these are homogeneous and were known since 1968 \cite{Gray:Wolf}: the round $6$--sphere endowed with the non-integrable almost complex structure induced by octonionic multiplication on $\R^7\simeq\Imag\mathbb{O}$ and the $3$--symmetric spaces $\CP^3=Sp(2)/U(1)\times Sp(1)$, $S^3\times S^3=SU(2)^3/SU(2)$ and $\F_3=SU(3)/T^2$. Recently two inhomogeneous nearly K\"ahler structures on $S^6$ and $S^3\times S^3$ were found in \cite{Foscolo:Haskins}. Finite quotients of the homogeneous nearly K\"ahler manifolds have also been studied \cite{Cortes}. This scarcity of examples should be contrasted with the infinitely many known examples of manifolds with real Killing spinors in other dimensions: Sasaki--Einstein, $3$--Sasaki and nearly parallel \gtwo  manifolds (the cross-sections of Calabi--Yau, hyperk\"ahler and \spin{7} cones, respectively).

In this paper we study the deformation theory of nearly K\"ahler manifolds. In \cite{Moroianu:Nagy:Semmelmann} Moroianu--Nagy--Semmelmann studied infinitesimal deformations of nearly K\"ahler structures, identifying them with an eigenspace of the Laplacian acting on coclosed primitive $(1,1)$--forms. The question of whether nearly K\"ahler $6$--manifolds have smooth, unobstructed deformations was however left open. Because of the scarcity of examples it would be very interesting to understand whether it is possible to obtain new nearly K\"ahler manifolds by deforming known ones. Understanding whether nearly K\"ahler deformations are in general obstructed is also important for applications to the theory of \gtwo \emph{conifolds} (asymptotically conical and conically singular \gtwo manifolds) developed by Karigiannis--Lotay \cite{Karigiannis:Lotay}. Finally, possible constructions of new examples of nearly K\"ahler manifolds based on singular perturbation methods require as a preliminary step the study of the deformation theory of nearly K\"ahler manifolds (and its extension to certain singular nearly K\"ahler spaces).

It is instructive to recall known results about the deformation theory of manifolds with real Killing spinors in other dimensions. Continuous families of Sasaki--Einstein structures are certainly known, \eg the regular Sasaki--Einstein $5$--manifolds obtained from del Pezzo surfaces of degree $d\leq 4$ via the Calabi ansatz have non-trivial moduli. However in general Sasaki--Einstein manifolds have obstructed deformations (\cf \cite{Tipler:vanCoevering} for the relation between integrability of infinitesimal deformations and $K$--stability in the more general context of constant scalar curvature Sasaki metrics). By a result of Pedersen--Poon \cite{Pedersen:Poon} $3$--Sasaki manifolds are rigid. In \cite{Alexandrov:Semmelmann} Alexandrov--Semmelmann study infinitesimal deformations of nearly parallel \gtwo structures. As in the nearly K\"ahler case these are identified with a certain subspace of an eigenspace of the Laplacian acting on $3$--forms. It is unclear whether infinitesimal deformations of nearly parallel \gtwo manifolds are unobstructed in general.

Given what is known about deformations of Sasaki--Einstein manifolds, it is not surprising that nearly K\"ahler $6$--manifolds have obstructed deformations in general, as we show in this paper.

In \cite[Theorem 6.12]{Koiso} Koiso showed that infinitesimal deformations of Einstein metrics are in general obstructed. He exhibited Einstein symmetric spaces with non-trivial infinitesimal Einstein deformations which cannot be integrated to second order. We will follow a similar strategy.

In \cite{Moroianu:Semmelmann} Moroianu--Semmelmann calculated the space of infinitesimal deformations of the homogeneous nearly K\"ahler manifolds. They found that $\CP^3$ and $S^3\times S^3$ are rigid while the flag manifold $\F_3$ has an $8$--dimensional space of infinitesimal deformations \cite[Corollary 6.1]{Moroianu:Semmelmann}. (The case of the round $6$--sphere is special, since there are more than two Killing spinors in this case. The space of nearly K\"ahler structures compatible with the round metric is an $\R\PP^7$--bundle over $S^6$. Since the round metric does not admit any Einstein deformation, there are no infinitesimal nearly K\"ahler deformations other than the ones coming from this family.)

In this paper we address the question of deciding whether the homogeneous nearly K\"ahler structure on the flag manifold admits genuine nearly K\"ahler deformations.

\begin{theorem*}\label{thm:Main:1}
The infinitesimal nearly K\"ahler deformations of the flag manifold $\F_3$ are all obstructed.
\end{theorem*}

The proof of the Theorem (and the paper itself) is divided into two distinct parts. In a first step we will obtain a deeper understanding of the deformation theory of nearly K\"ahler structures in general beyond the infinitesimal level considered in \cite{Moroianu:Nagy:Semmelmann}. The main tool is a certain Dirac-type operator on nearly K\"ahler manifolds, which appears as a certain combination of differential, codifferential and type decomposition acting on differential forms. The use of Dirac operators as tools in the deformation theory of manifolds with special geometric structures is not new. N\"ordstrom \cite[Chapter 3]{Nordstrom:Thesis} used Dirac operators to streamline certain steps in the deformation theory of manifolds with special holonomy. More precisely, mapping properties of Dirac operators are used to establish slice theorems for the action of the diffeomorphism group. This approach has turned out to be particularly useful in the deformation theory of \emph{non-compact} manifolds with special holonomy, in particular asymptotically cylindrical manifolds \cite{Nordstrom:ACyl} and conifolds \cite{Karigiannis:Lotay}.

Besides the application to the deformation theory of nearly K\"ahler structures, we will show that Dirac-type operators can also be used to obtain interesting results about Hodge theory on nearly K\"ahler manifolds. In particular, we will give an elementary proof of a result of Verbitsky \cite[Theorem 6.2 and Remark 6.4]{Verbitsky} on the type decomposition of harmonic forms on nearly K\"ahler manifolds.

A second important ingredient in our treatment of the deformation theory of nearly K\"ahler manifolds is Hitchin's interpretation of nearly K\"ahler structures as (constrained) critical points of a Hamiltonian function on the infinite dimensional symplectic vector space $\Omega^3_{exact}\times \Omega^4_{exact}$. This description will allow us to interpret the nearly K\"ahler equations as the vanishing of a smooth map $\Phi\co \Omega^3_{exact}\times \Omega^4_{exact}\times\Omega^1\ra \Omega^3_{exact}\times \Omega^4_{exact}$, rather than as equations on the space of \sunitary{3} structures. The main advantage of this approach is to introduce additional free parameters that can be used to reduce the cokernel of the linearisation $D\Phi$ and pin down exactly where possible obstructions to integrate infinitesimal deformations could lie.

The general deformation theory of nearly K\"ahler structures is then applied to the specific case of the flag manifold $\F_3$. The framework introduced in the first part of the paper will make possible the explicit calculation of the non-vanishing obstructions to integrate the infinitesimal deformations of the homogeneous nearly K\"ahler structure on $\F_3$ to second order. 

As we have already mentioned, Alexandrov--Semmelmann \cite{Alexandrov:Semmelmann} studied infinitesimal deformations of nearly parallel \gtwo manifolds. In particular, by \cite[\S 8]{Alexandrov:Semmelmann} the normal homogeneous nearly parallel \gtwo manifolds are all rigid except for the Aloff--Wallach manifold $SU(3)/U(1)$, which admits an $8$--dimensional space of infinitesimal deformations isomorphic to $\Lie{su}_3$. It is likely that the methods of this paper could also be used to analyse the integrability of these infinitesimal deformations.

\vspace{0.2cm}

The paper is organised as follows. In Section 2 we collect various preliminary results about $6$--manifolds endowed with an \sunitary{3} structure. The notion of a \emph{stable} form introduced by Hitchin in \cite{Hitchin} will be central in our exposition. These results are known and we collect them here simply for the convenience of the reader. In Section 3 we study algebraic compatibilities on nearly K\"ahler manifolds between differential and codifferential and the decomposition of the space of differential forms into types corresponding to irreducible representations of \sunitary{3}. We will then introduce the Dirac-type operator mentioned above and study its mapping properties. As a first application, we will derive results about the Hodge theory of nearly K\"ahler manifolds. Section 4 discusses the deformation problem of nearly K\"ahler manifolds. The Dirac-type operator is used to define a slice for the action of the diffeomorphism group while Hitchin's interpretation of nearly K\"ahler structures as (constrained) critical points allows to re-write the nearly K\"ahler equations as the vanishing of a certain non-linear map $\Phi\co \Omega^3_{exact}\times \Omega^4_{exact}\times\Omega^1\ra \Omega^3_{exact}\times \Omega^4_{exact}$. We study the linearisation $D\Phi$ at a nearly K\"ahler structure and identify its cokernel and therefore possible obstructions to integrate infinitesimal deformations to genuine nearly K\"ahler deformations. Deformations of the homogeneous nearly K\"ahler structure on the flag manifold are finally studied in Section 5. We introduce the notion of second order deformations, recall the result of  Moroianu--Semmelmann \cite[Corollary 6.1]{Moroianu:Semmelmann} on the existence of an $8$--dimensional space of infinitesimal nearly K\"ahler deformations and, via explicit calculations and representation theoretic considerations, we show that these are all obstructed to second order.

\vspace{-0.2cm}
  
\subsection*{Acknowledgements} The author wishes to thank Johannes Nordstr\"om and Uwe Semmelmann for interesting conversations related to this paper and Uwe Semmelmann for comments on an earlier draft of this note. Part of this work was carried out while the author was visiting Leibniz Universit\"at Hannover in the Fall 2015; he wishes to thank the Riemann Center for Geometry and Physics for support and the Differential Geometry Group for hospitality. 

\section{\sunitary{3} structures on $6$--manifolds}

In this preliminary section we collect various known facts about $SU(3)$ structures on $6$--manifolds. These results are well-known and we collect them here simply for the convenience of the reader.

\subsection{Stable forms}

Following Hitchin \cite{Hitchin}, the notion of a stable form will be central in our exposition.

\begin{definition}\label{def:Stable:Forms}
A differential form $\phi \in \Lambda^p(\R^n)^\ast$ is \emph{stable} if its $GL(n,\R)$--orbit in $\Lambda^p(\R^n)^\ast$ is open.
\end{definition}

In dimension $6$, there are only three possibilities for stable forms \cite[\S 2]{Hitchin}:
\begin{enumerate}
\item a stable $2$--form $\omega$ (a non-degenerate $2$--form) with open orbit isomorphic to $GL(6,\R)/Sp(6,\R)$;
\item a stable $4$--form $\sigma$, with stabiliser $Sp(6,\R)$;
\item a stable $3$--form $\rho$, with stabiliser $SL(3,\C)$.
\end{enumerate}

Note that in all three cases the stabiliser is in fact contained in $SL(6,\R)$ and therefore to each stable form one can associate a volume form $\dvol$. For example one can define $\dvol (\omega) = \frac{1}{n!}\omega^n$ for a stable $2$--form $\omega$ in dimension $2n$. Using the homogeneous behaviour of the map $\dvol$, for every stable $p$--form $\phi$ Hitchin defines its \emph{dual} $\hat{\phi}$, a $(n-p)$--form such that $\dvol (\phi)$ is proportional to $\phi \wedge\hat{\phi}$. In dimension $6$ we have:
\begin{enumerate}
\item for a stable $2$--form $\omega$, $\hat{\omega}=\frac{1}{2}\omega^{2}$;
\item for a stable $4$--form $\sigma$, $\hat{\sigma}$ is the unique non-degenerate $2$--form such that $\sigma = \frac{1}{2}\hat{\sigma}^2$;
\item for a stable $3$--form $\rho$, $\hat{\rho}$ is the unique $3$--form such that $\rho + i \hat{\rho}$ is  a nowhere vanishing complex volume form.
\end{enumerate}
In particular, in dimension $6$ the real part of a complex volume form $\Omega$ uniquely determines $\Omega$ itself. Moreover, since its stabiliser is $SL(3,\C)$, a stable $3$--form $\Real\Omega$ defines an almost complex structure $J$: a $1$--form $\alpha$ is of type $(1,0)$ if and only if $\alpha\wedge\Omega=0$.

\begin{definition}\label{def:SU(3):structure}
An $SU(3)$--structure on a $6$--manifold $M$ is a pair of smooth differential forms $(\omega,\Real\Omega)$ such that
\begin{enumerate}
\item $\omega$ is a stable $2$--form with Hitchin's dual $\frac{1}{2}\omega^2$;
\item $\Real\Omega$ is a stable $3$--form with Hitchin's dual $\Imag\Omega$;
\item the following algebraic constraints are satisfied:
\begin{equation}\label{eq:SU(3):structure:Constraints}
\omega\wedge\Real\Omega=0, \qquad \tfrac{1}{6}\omega^3 = \tfrac{1}{4}\Real\Omega\wedge\Imag\Omega.
\end{equation}
\end{enumerate}
\end{definition}
The two constraints guarantee that the stabiliser of the pair $(\omega,\Real\Omega)$ is exactly $SU(3)=Sp(6,\R)\cap SL(3,\C)$. We could of course define an $SU(3)$--structure as the choice of a pair of stable differential forms $\sigma\in\Omega^4(M)$ and $\Real\Omega\in \Omega^3 (M)$ satisfying \eqref{eq:SU(3):structure:Constraints} with $\omega=\hat{\sigma}$.

On $\R^6 \simeq \C^3$ with holomorphic coordinates $(z_1,z_2,z_3)$ we define the standard parallel $SU(3)$--structure $(\omega_0, \Real\Omega_0)$ by
\[
\omega_0 = \tfrac{i}{2} \left( dz_1 \wedge d\overline{z}_1 + dz_2 \wedge d\overline{z}_2 + dz_3 \wedge d\overline{z}_3 \right), \qquad \Omega_0 = dz_1 \wedge dz_2 \wedge dz_3. 
\]
An $SU(3)$ structure on a $6$--manifold $M$ in the sense of Definition \ref{def:SU(3):structure} defines a reduction of the structure group of the tangent bundle of $M$ to $SU(3)$ by considering the sub-bundle of the frame-bundle of $M$ defined by $\{ u\co \R^6 \stackrel{\sim}{\longrightarrow} T_{x}M \, |\, u^\ast (\omega_x,\Real\Omega_x)=(\omega_0,\Real\Omega_0) \}$.

\begin{remark}\label{rmk:SU(3):spin}
Since $SU(3)\subset SU(4)\simeq Spin(6)$ is precisely the stabiliser of a non-zero vector in $\C^4$, we could also define an $SU(3)$ structure as the choice of a spin structure on $M$ together with a non-vanishing spinor.
\end{remark}

Note that every $SU(3)$--structure induces a Riemannian metric $g$ because $SU(3)\subset SO(6)$. Hence from now on we will identify without further notice $TM$ and $T^\ast M$ using the metric $g$ (and $\R^6$ with $(\R^6)^\ast$ using the flat metric $g_0$).

\begin{remark}\label{rmk:T:Tast}
In particular, when we write $JX$ and think of it as a $1$--form we really mean $(JX)^\flat$. The fact that $(JX)^\flat = -JX^\flat$ might cause some confusion at times. We will instead distinguish between differential $df$ and gradient $\nabla f$ of a function $f$.
\end{remark}

\subsection{Decomposition of the space of differential forms}

The decomposition into irreducible representations of the $SU(3)$--representation $\Lambda^\ast\R^6$ is well-known. This is usually stated after complexification in terms of the $(p,q)$--type decomposition induced by the standard complex structure $J_0$ and in terms of primitive forms. We will stick with real representations and find more convenient to use the uniform notation $\Lambda^k_\ell$ for an irreducible component of $\Lambda^k\R^6$ of dimension $\ell$.

\begin{lemma}\label{lem:decomposition:forms}
We have the following orthogonal decompositions into irreducible $SU(3)$ representations:
\[
\Lambda ^{2}\R^6= \Lambda^2_1 \oplus \Lambda^2_6 \oplus \Lambda^2_8,
\]
where $\Lambda^2_1 = \R\,\omega$, $\Lambda^2_6 = \{ X\lrcorner\Real\Omega \, |\, X\in\R^6\}$ and $\Lambda^2_8$ is the space of primitive forms of type $(1,1)$.
\[
\Lambda ^{3}\R^6= \Lambda^3 _6 \oplus \Lambda^3 _{1\oplus 1} \oplus \Lambda^3_{12},
\]
where $\Lambda^3_6 = \{ X\wedge \omega \, |\, X\in\R^6\}$, $\Lambda^3 _{1\oplus 1} = \R \Real\Omega \oplus \R \Imag\Omega$ and $\Lambda^3_{12}$ is the space of primitive forms of type $(1,2)+(2,1)$, \ie $\Lambda^3_{12} = \{ S_{\ast}\Real{\Omega} \, | \, S \in Sym^{2}(\R ^{6}), SJ+JS=0 \}$.
\end{lemma}
\begin{remark}
Here an endomorphism $S\in \text{End}(\R^6)$ acts on a differential $p$--form $\phi$ by $S_\ast \phi (X_1,\dots, X_p) = -\sum_{j=1}^p{\phi(X_1,\dots, SX_j,\dots,X_p)}$.
\end{remark}
For a $6$--manifold $M$ endowed with an $SU(3)$--structure $(\omega,\Omega)$ we denote by $\Omega^k_\ell (M)$ the space of smooth sections of the bundle over $M$ with fibre $\Lambda^k_\ell$.

By definition $\Lambda^2_6$ is isomorphic to $\Lambda^1$ via the metric and the contraction with $\Real\Omega$. The adjoint of this map with respect to the flat metric $g_0$ will be denoted by $\alpha\co \Lambda^2_6\ra\Lambda^1$. Note that for a form $\eta\in\Lambda^2_6$ we have
\begin{equation}\label{eq:alpha}
\eta=\tfrac{1}{2}\alpha(\eta)\lrcorner\Real\Omega.
\end{equation}

The following identities follow from \cite[Equations (12), (17), (18) and (19)]{Moroianu:Nagy:Semmelmann}.

\begin{lemma}\label{lem:Hodge:star}
In the decomposition of Lemma \ref{lem:decomposition:forms} the Hodge--$\ast$ operator is given by:
\begin{enumerate}
\item $\ast\omega = \tfrac{1}{2}\omega^2$;
\item $\ast (X\lrcorner\Real\Omega) = -JX\wedge\Real\Omega = X\wedge\Imag\Omega$;
\item $\ast (\eta _{0} \wedge \omega)=-\eta _{0}$ for all $\eta _{0} \in \Omega ^2_8$;
\item $\ast (X \wedge \omega) = \frac{1}{2} X \lrcorner \omega ^{2}=JX \wedge \omega$;
\item $\ast \Real{\Omega}=\Imag{\Omega}$ and $\ast\Imag\Omega = -\Real\Omega$;
\item $\ast (S_{\ast}\Real{\Omega})=-S_{\ast}\Imag{\Omega}=(JS)_{\ast}\Real{\Omega}$;
\end{enumerate}
\end{lemma}

We use Lemma \ref{lem:Hodge:star} to deduce useful identities and characterisations of the different types of forms.

\begin{lemma}\label{lem:identities:2forms}
If $\eta = \eta _{0}+\lambda \omega + X \lrcorner\Real{\Omega} \in \Omega ^{2}$ with $\eta_0\in\Omega^2_8$, then the following holds:
\begin{enumerate}
\item $\ast (\eta \wedge \omega) = -\eta _{0}+2\lambda \omega + X \lrcorner \Real{\Omega}$;
\item $\ast( \eta \wedge \eta \wedge \omega ) = -|\eta _{0}|^{2}+6\lambda ^{2}+2|X|^{2}$;
\item $\ast (\eta \wedge \Real{\Omega}) = 2JX$ and $\ast (\eta \wedge \Imag{\Omega})=-2X$;
\item $\ast (\eta \wedge \omega ^{2}) = 6\lambda$.
\end{enumerate}
In particular, $\eta \in \Omega ^2_8$ iff $\eta \wedge \omega^{2} =0=\eta \wedge \Real{\Omega}$ iff $\eta \wedge \omega =-\ast \eta$.
\proof
(i) follows immediately from Lemma \ref{lem:Hodge:star} (i)--(iii). In turn, (i) implies (ii) by the identity
\[
\eta \wedge \eta \wedge \omega = \eta \wedge \ast ^{2}(\eta \wedge \omega),
\]
the fact that the decomposition of Lemma \ref{lem:decomposition:forms} is orthogonal and $|\omega|^{2}=3$, $\langle X \lrcorner\Real{\Omega}, X \lrcorner\Real{\Omega} \rangle =2|X|^{2}$. The identities (iii) follows immediately from \cite[Equations (3) and (4)]{Moroianu:Nagy:Semmelmann} and (iv) from $2\ast\omega = \omega ^{2}$ and $|\omega|^{2}=3$.
\endproof
\end{lemma}

We have similar identities on $3$--forms, with analogous proof.

\begin{lemma}\label{lem:identities:3forms}
If $\sigma = X \wedge \omega + \lambda \Real{\Omega} + \mu \Imag{\Omega}+ S_{\ast}\Real{\Omega} \in \Omega ^{3}$, then the following holds:
\begin{enumerate}
\item $\ast (\sigma \wedge \omega) = 2JX$;
\item $\ast (\sigma \wedge \Real{\Omega}) = -4\mu$;
\item $\ast (\sigma \wedge \Imag{\Omega})=4\lambda$.
\end{enumerate}
In particular, $\rho \in \Omega ^3_{12}$ iff $\rho \wedge \omega =0=\rho \wedge \Omega$.
\end{lemma}

Finally, it will be useful to have an explicit formula for the linearisation of Hitchin's duality map for stable forms in terms of the Hodge--$\ast$ and the decomposition of forms into types.

\begin{prop}\label{prop:Linearisation:Hitchin:dual}
Given an $SU(3)$ structure $(\omega,\Real\Omega)$ on $M^6$ let $\sigma\in\Omega^4(M)$ and $\rho\in\Omega^3 (M)$ be forms with small enough $C^0$--norm so that $\frac{1}{2}\omega^2+\sigma$ and $\Real\Omega+\rho$ are still stable forms. Decomposing into types we write $\sigma = \sigma_1 + \sigma _6 + \sigma_8$ and $\rho = \rho _6 + \rho _{1\oplus 1} + \rho _{12}$.
\begin{enumerate}
\item The image $\hat{\sigma}$ of $\sigma$ under the linearisation of Hitchin's duality map at $\tfrac{1}{2}\omega^2$ is
\[
\hat{\sigma} = \tfrac{1}{2}\ast\sigma_1 + \ast \sigma_6 -\ast\sigma_8.
\]
\item The image $\hat{\rho}$ of $\rho$ under the linearisation of Hitchin's duality map at $\Real\Omega$ is
\[
\hat{\rho} = \ast (\rho_6 + \rho _{1\oplus 1}) -\ast \rho_{12}. 
\]
\end{enumerate}
\proof
In order to prove the first statement, observe that $\hat{\sigma}$ is the unique $2$--form such that $\sigma = \omega\wedge\hat{\sigma}$. Apply $\ast$ to this identity and use Lemma \ref{lem:identities:2forms}.(i). The formula for $\hat{\rho}$ follows from \cite[Lemma 3.3]{Moroianu:Nagy:Semmelmann} and the last three identities in Lemma \ref{lem:Hodge:star}.
\endproof
\end{prop}

\subsection{Nearly K\"ahler manifolds}

Given a subgroup $G$ of $SO(n)$, we define a $G$--structure on a Riemannian manifold $(M^n,g)$ as a sub-bundle $\mathcal{P}$ of the orthogonal frame bundle of $M$ with structure group $G$. The \emph{intrinsic torsion} of $\mathcal{P}$ is a measure of how much $\mathcal{P}$ is far from being parallel with respect to the Levi--Civita connection $\nabla$ of $(M,g)$. More precisely, restricting $\nabla$ to $\mathcal{P}$ yields a $\Lie{so}(n)$--valued $1$--form $\theta$ on $\mathcal{P}$. Choose a complement $\Lie{m}$ of the Lie algebra of $G$ in $\Lie{so}(n)$. Projection of $\theta$ onto $\Lie{m}$ yields a $1$--form $T$ on $M$ with values in the bundle $\mathcal{P}\times_G \Lie{m}$. This is the (intrinsic) torsion of the $G$--structure $\mathcal{P}$.

For an $SU(3)$--structure on a $6$--manifold $M$ one can check that $\nabla (\omega,\Omega) = T_\ast (\omega,\Omega)$ where $T$ acts on differential forms via the representation of $\Lie{m}\subset\Lie{so}(6)$ on $\Lambda^\ast(\R^6)$. It turns out that $T$ itself is uniquely recovered by knowledge of the anti-symmetric part of $T_\ast (\omega,\Omega)$, \ie the knowledge of $d\omega$ and $d\Omega$.

\begin{prop}[{\cite[Theorem 2.9]{Bedulli:Vezzoni:SU(3)}}]\label{prop:Torsion:SU(3):structures}
Let $(\omega,\Omega)$ be an $SU(3)$ structure. Then there exists functions $w_1,\hat{w}_1 \in \Omega^0$, $w_2,\hat{w}_2\in\Omega^2_8$, $w_3\in\Omega^3_{12}$ and vector fields $w_4,w_5$ on $M$ such that
\[
\begin{gathered}
d\omega = 3w_1 \Real\Omega + 3\hat{w}_1 \Imag\Omega + w_3 + w_4 \wedge\omega,\\
d\Real\Omega = 2\hat{w}_1\omega^2 + w_5\wedge\Real\Omega + w_2 \wedge\omega,\\
d\Imag\Omega=-2w_1\omega^2 - Jw_5 \wedge\Real\Omega + \hat{w}_2\wedge\omega.
\end{gathered}
\]
\end{prop}
Note that the different sign in front of $Jw_5\wedge\Real\Omega$ in the formula for $d\Imag\Omega$ with respect to the formula in \cite[Theorem 2.9]{Bedulli:Vezzoni:SU(3)} is consistent with Remark \ref{rmk:T:Tast}.

\begin{definition}\label{def:NK}
An $SU(3)$ structure on a $6$--manifold $M$ is called a \emph{nearly K\"ahler} structure if $\hat{w}_1,w_2,\hat{w}_2,w_3,w_4,w_5$ all vanish and $w_1 =1$. In other words a nearly K\"ahler structure satisfies
\begin{equation}\label{eq:NK}
d\omega = 3\Real\Omega, \qquad d\Imag\Omega = -2\omega^2.
\end{equation}
\end{definition}

As remarked in the Introduction, \eqref{eq:NK} are equivalent to the requirement that
\[
\varphi = r^2 dr\wedge\omega + r^3 \Real\Omega
\]
is a closed and coclosed ``conical'' $G_2$ structure on the cone $C(M)=\R^+\times M$. Thus the cone $C(M)$ has holonomy contained in $G_2$ and in particular is Ricci-flat. As a consequence, nearly K\"ahler manifolds are Einstein with positive Einstein constant normalised so that $\text{Scal}=30$. In particular every complete nearly K\"ahler manifold is compact with finite fundamental group.

In Remark \ref{rmk:SU(3):spin} we observed that every $6$--manifold with an $SU(3)$ structure is spin and endowed with a unit spinor $\psi$. The nearly K\"ahler equations \eqref{eq:NK} have an equivalent interpretation as a first order differential equation on $\psi$ \cite[Theorem 2]{Bar}. Indeed, every $G_2$ manifold admits a parallel spinor. Restricting the parallel spinor on the cone $C(M)$ to the cross-section $M$ induces a \emph{real Killing spinor} $\psi$ on every nearly K\"ahler manifold, \ie (possibly after changing orientation) a unit spinor $\psi$ such that
\begin{equation}\label{eq:Killing:spinor}
\nabla _X \psi = \tfrac{1}{2}X\cdot\psi
\end{equation}
for every vector field $X$. Clifford multiplication by the volume form $\vol$ (\ie the complex structure on the spinor bundle) yields a second Killing spinor $\vol\cdot\psi$ satisfying (recall that $X\cdot\vol\cdot\psi=-\vol\cdot X\cdot\psi$)
\[
\nabla _X (\vol\cdot\psi) =-\tfrac{1}{2}X\cdot(\vol\cdot\psi).
\]

In the rest of the paper $(M,\omega,\Omega)$ will denote a complete (hence compact) nearly K\"ahler $6$--manifold and $\psi$ will denote the real Killing spinor on $M$ satisfying \eqref{eq:Killing:spinor}.

\section{Hodge theory on nearly K\"ahler manifolds}

The main goal of this section is to introduce a Dirac-type operator $D$ on a nearly K\"ahler manifold $M$ and study its mapping properties. $D$ differs from the standard Dirac operator by a lower order term. This operator arises as a certain composition of differential, codifferential and type decomposition on differential forms and will play a central role in the study of the deformation theory of nearly K\"ahler manifolds. Furthermore, we find that the mapping properties of the operator $D$ have interesting consequences on the Hodge theory of nearly K\"ahler manifolds. In particular, we will give an elementary proof of results of Verbitsky \cite{Verbitsky} on the type decomposition of harmonic forms on nearly K\"ahler manifolds.

\subsection{Differential and co-differential on nearly K\"ahler manifolds}

Before giving the definition of the Dirac operator $D$ we need to study how the exterior differential $d\co \Omega^k(M) \ra \Omega^{k+1}(M)$ and its adjoint behave with respect to the decomposition of $\Omega^\ast(M)$ introduced in Lemma \ref{lem:decomposition:forms}. We need an additional piece of notation and few simple observations.

Let $\Lambda\co \Omega^k(M)\ra\Omega^{k-2}(M)$ be the point-wise adjoint of wedging with $\omega$ and note that $\Lambda (X\wedge\omega)=2X$ for every vector field $X$. Recall also that $\alpha\co \Omega^3(M) \ra \Omega^1(M)$ is the point-wise adjoint of the map $X\mapsto X\lrcorner\Real\Omega$. Furthermore we have the following identities:
\begin{equation}\label{eq:Hodge:star:contraction}
\ast (X\wedge\sigma)=(-1)^p X\lrcorner\ast\sigma
\end{equation}
for every $p$--form $\sigma$ and $X\in \mathcal{X}(M)$;
\begin{equation}\label{eq:identities:1:forms:a}
(X\lrcorner\Real\Omega)\wedge\omega = -JX\wedge\Real\Omega, 
\end{equation}
which follows immediately by contracting $\Real\Omega\wedge\omega=0$ by $X$.
\begin{equation}\label{eq:identities:1:forms:b}
(X\lrcorner\Real\Omega)\wedge\Real\Omega = X\wedge\omega^2, \qquad (X\lrcorner\Real\Omega)\wedge\Imag\Omega = JX\wedge\omega^2
\end{equation}
which are \cite[Equations (3) and (4)]{Moroianu:Nagy:Semmelmann}.
\begin{equation}\label{eq:identities:1:forms:c}
\ast Y = \tfrac{1}{2}JY\wedge\omega^2
\end{equation}
which is a consequence of \eqref{eq:Hodge:star:contraction}. As a consequence,
\begin{equation}\label{eq:Dirac:prelim:identities}
d^\ast Y = -\ast (dJY\wedge\tfrac{1}{2}\omega^2).
\end{equation}

\begin{prop}\label{prop:differential:2forms}
For every $f\in C^\infty(M)$, $\eta \in \Omega^2_8 (M)$ and $X\in\mathcal{X}(M)$ we have
\begin{enumerate}
\item $d(f\omega)=df\wedge\omega+ 3f\Real\Omega$;
\item $d^\ast (f\omega) = J df$;
\item $d\eta = \tfrac{1}{2}Jd^{\ast}\eta \wedge \omega + \rho$ for some $\rho \in \Omega ^3_{12}$;
\item $d(X \lrcorner \Real{\Omega}) = -\left( \tfrac{1}{2}\alpha (dJX) + 3X \right) \wedge \omega -\tfrac{1}{2}(d^{\ast}X) \Real{\Omega} -\tfrac{1}{2}d^{\ast}(JX) \Imag{\Omega} + \rho'$ for some $\rho' \in \Omega ^3_{12}$;
\item $d^\ast(X\lrcorner\Real{\Omega})=J\alpha (dJX)$.
\end{enumerate}
\proof
(i) follows immediately from \eqref{eq:NK}.

Since $\omega^2$ is closed, $d^\ast(f\omega) = -\ast d(\tfrac{1}{2}f\omega^2)=-\ast(df\wedge \tfrac{1}{2}\omega^2)$ and \eqref{eq:identities:1:forms:c} completes the proof of (ii).

In order to prove (iii), differentiate $\eta \wedge \Omega =0$ and use the fact that $\eta \wedge d\Omega =-2i\eta\wedge\omega^2 =0$. Lemma \ref{lem:identities:3forms} then implies that $d\eta$ has zero component in the complex line spanned by $\Omega$. Similarly, differentiating the equality $\eta \wedge \omega + \ast \eta=0$ and using Lemma \ref{lem:identities:3forms}.(i) yields $\Lambda d\eta=Jd^{\ast}\eta$.

The identity (iv) is proved in a similar way. Consider $d(X \lrcorner \Real{\Omega}) \wedge \omega$: by \eqref{eq:identities:1:forms:a} and \eqref{eq:identities:1:forms:b}
\[
d(X\lrcorner\Real\Omega)\wedge\omega=-3X\wedge\omega^2 -dJX\wedge\Real\Omega.
\]
Lemma \ref{lem:identities:2forms}.(iii) and \eqref{eq:identities:1:forms:c} now imply that the $\Omega^3_6$--component of $d(X\lrcorner\Real\Omega)$ is
\[
-\left( \tfrac{1}{2}\alpha (dJX) + 3X \right) \wedge \omega
\]
as claimed. On the other hand, by \eqref{eq:identities:1:forms:b} and the fact that $\omega^2$ is closed
\[
d(X\lrcorner\Real\Omega)\wedge\Real\Omega = dX\wedge\omega^2, \qquad d(X\lrcorner\Real\Omega)\wedge\Imag\Omega = dJX\wedge\omega^2.
\]
(iv) now follows from \eqref{eq:Dirac:prelim:identities} and Lemma \ref{lem:identities:2forms}.(iv).

In order to prove (v) observe that $X\lrcorner\Real\Omega = JX\lrcorner\Imag\Omega$ and therefore by \eqref{eq:Hodge:star:contraction},
\[
d^\ast (X\lrcorner\Real\Omega)=\ast d(JX\wedge\Real\Omega)=\ast (dJX\wedge\Real\Omega).
\]
Lemma \eqref{lem:identities:2forms}.(iii) concludes the proof.
\endproof
\end{prop}

We actually need a bit more: the $\Omega^3_{12}$ component of $d(X\lrcorner\Real\Omega)$ is equal to the $\Omega^3_{12}$ component of $d^\ast(X\wedge\Real\Omega)$. In particular, if $d(X\lrcorner\Real\Omega)\in\Omega^3_{12}$ then $d(X\lrcorner\Real\Omega)=0$.

\begin{lemma}\label{lem:Lie:derivative}
For every vector field $X \in \mathcal{X}(M)$ we have
\[
d(X\lrcorner\Real\Omega)=\left( J\alpha (dX)+2X\right) \wedge\omega - d^\ast X\Real\Omega -d^\ast(JX)\Imag\Omega + d^\ast(X\wedge\Real\Omega).
\]
\proof
If $\rho=\mathcal{L}_X\Real\Omega$ then the linearisation of Hitchin's duality map $\hat{\rho}$ must be
\[
\hat{\rho}=\mathcal{L}_X\Imag\Omega=d(X\lrcorner\Imag\Omega)-4JX\wedge\omega=-d(JX\lrcorner\Real\Omega)-4JX\wedge\omega.
\]
On the other hand, by Proposition \ref{prop:differential:2forms}.(iv)
\[
\rho = d(X\lrcorner\Real\Omega)=-\left( \tfrac{1}{2}\alpha (dJX) + 3X \right) \wedge \omega -\tfrac{1}{2}d^{\ast}X \Real{\Omega} -\tfrac{1}{2}d^{\ast}(JX) \Imag{\Omega} + \rho _{0}
\]
and Proposition \ref{prop:Linearisation:Hitchin:dual} now implies
\[\begin{gathered}
\hat{\rho}=-\left( \tfrac{1}{2}J\alpha (dJX) + 3JX \right) \wedge \omega +\tfrac{1}{2}d^{\ast}(JX) \Real{\Omega} -\tfrac{1}{2}d^{\ast}X \Imag{\Omega}  -\ast\rho _{0}\\
=-\left( J\alpha (dJX) + 6JX \right) \wedge \omega +d^{\ast}(JX) \Real{\Omega} -d^{\ast}X \Imag{\Omega}  -\ast d(X\lrcorner\Real\Omega).
\end{gathered}\]

Comparing the two expressions for $\hat{\rho}$ we conclude that
\[
d(JX\lrcorner\Real\Omega)=\left( J\alpha (dJX) + 2JX \right) \wedge \omega -d^{\ast}(JX) \Real{\Omega} +d^{\ast}X \Imag{\Omega}  +\ast d(X\lrcorner\Real\Omega).
\]
Up to changing $X$ into $JX$, the Lemma is now proved since by \eqref{eq:Hodge:star:contraction}
\[
d^\ast (JX\wedge\Real\Omega)=\ast d(X\lrcorner\Real\Omega).\qedhere
\]
\end{lemma}

\begin{remark}\label{rmk:Lie:derivatives}
Equating the two different ways of writing $\hat{\rho}$ also yields the identity
\begin{equation}\label{eq:dX:dJX}
\alpha(dX)=J\alpha(dJX)+4JX
\end{equation}
of \cite[Lemma 3.2]{Moroianu:Semmelmann}. Note also that integrating by parts the identity of Lemma \ref{lem:Lie:derivative} with a compactly supported $X$ yields
\begin{equation}\label{eq:alpha:primitive:(2,1)}
\alpha(d^\ast\rho)=-J \alpha(\ast d\rho)
\end{equation}
for every $\rho\in\Omega^3_{12}(M)$. In particular, if $\rho\in\Omega^3_{12}$ is coclosed then $d\rho\in \Omega^4_8$. Indeed, since $\rho\wedge\omega=0=\rho\wedge\Real\Omega$ by Lemma \ref{lem:identities:3forms} we also know that $d\rho$ has not component in $\Omega^4_1$. 
\end{remark}

\subsection{The Dirac operator on nearly K\"ahler manifolds}

Every $6$--manifold $M$ with an $SU(3)$ structure is spin and it is endowed with a unit spinor $\psi$. As an $SU(3)$ representation the real spinor bundle $\slashed{S}(M)$ is isomorphic to $\Lambda^0\oplus\Lambda^0\oplus\Lambda^1$, where the isomorphism is
\[
(f,g,X)\longmapsto f\psi + g\vol\cdot\psi + X\cdot\psi.
\]
For comparison with the Dirac-type operator we will define in the next subsection, we want now to derive a formula for the Dirac operator $\slashed{D}$ on a nearly K\"ahler manifold in terms of this isomorphism.

Recall that on a nearly K\"ahler manifold we can assume the defining unit spinor $\psi$ satisfies \eqref{eq:Killing:spinor}. In particular, $\slashed{D}\psi = -3\psi$ and $\slashed{D}(\vol\cdot\psi)=3(\vol\cdot\psi)$. Thus
\[
\slashed{D}\left( f\psi + g\vol\cdot\psi\right) = -3f\psi + 3g\vol\cdot\psi +\left( \nabla f -J\nabla g \right)\cdot \psi,
\]
since the almost complex structure $J$ satisfies
\begin{equation}\label{eq:J:spinors}
JX\cdot\psi = \vol\cdot X\cdot\psi=-X\cdot\vol\cdot\psi.
\end{equation}

On the other hand,
\[
\slashed{D}(X\cdot\psi) = \sum_{i=1}^6{e_i\cdot\nabla_{e_i}X\cdot\psi}-X\cdot\psi -X\cdot\slashed{D}\psi =dX\cdot\psi + (d^\ast X)\, \psi + 2X\cdot\psi.
\]

Now decompose $dX$ into types: $dX=\tfrac{1}{3}d^\ast (JX)\,\omega + \tfrac{1}{2}\alpha (dX)\lrcorner\Real\Omega + \pi_8(dX)$. We have to determine the action of $2$--forms on $\psi$.

\begin{lemma}\label{lem:Clifford:multiplication:2forms}
For any $2$--form $\eta=\lambda\omega + Y\lrcorner\Real\Omega + \eta_0$ with $\eta_0\in\Omega^2_8$ we have
\[
\eta\cdot\psi = 3\lambda \vol\cdot\psi + 2JY\cdot\psi.
\]
\proof
Forms of type $\Omega^2_8$ are exactly those that act trivially on $\psi_0$. On the other hand, writing $\omega = \sum_{i=1}^3 {e_i\wedge Je_i}$ for an $SU(3)$--adapted orthonormal co-frame $\{ e_1, Je_1, e_2, Je_2, e_3, Je_3\}$,
\[
\omega\cdot\psi = \sum_{i=1}^3{(e_i\wedge Je_i)\cdot\psi} = \sum_{i=1}^3{e_i\cdot Je_i\cdot\psi} = -\sum_{i=1}^3{e_i\cdot e_i\cdot\vol\cdot\psi}=3\vol\cdot\psi.
\]
Here we used \cite[Equation (1.3)]{Friedrich:al}, the fact that $e_i$ is orthogonal to $Je_i$ and \eqref{eq:J:spinors}.

In order to calculate the action of $\Omega^2_6$, observe that the intrinsic torsion of a nearly K\"ahler structure is $-\Imag\Omega$ \cite[p. 3]{Moroianu:Semmelmann} and that this acts as $4$ on $\psi$ and annihilates spinors of the form $X\cdot\psi$, \cf for example \cite[Lemma 2]{Charbonneau:Harland}. Thus, using \cite[Equations (1.3) and (1.4)]{Friedrich:al},
\[
(Y\lrcorner\Real\Omega)\cdot\psi = (JY\lrcorner\Imag\Omega)\cdot\psi = -\tfrac{1}{2}\left( JY\cdot \Imag\Omega + \Imag\Omega\cdot JY\right)\cdot\psi = 2JY\cdot\psi. \qedhere
\]
\end{lemma}

\begin{prop}\label{prop:Dirac:Riemannian}
For every $f,g\in C^\infty (M)$ and $X\in\mathcal{X}(M)$ we have
\[
\slashed{D}(f\psi + g\vol\cdot \psi + X\cdot\psi) = (d^\ast X-3f)\, \psi + (d^\ast JX+3g)\vol\cdot\psi + \left( \nabla f - J\nabla g -\alpha (dJX) - 2X\right)\cdot\psi.
\]
\proof
Use Lemma \ref{lem:Clifford:multiplication:2forms} and \eqref{eq:dX:dJX} to rewrite $J\alpha (dX)$ as $-\alpha (dJX) - 4X$.
\endproof
\end{prop}

\subsection{A Dirac-type operator}

Consider now the first order operator
\[
D\co \Omega^1\oplus\Omega^0\oplus\Omega^0\ra\Omega^1\oplus\Omega^0\oplus\Omega^0
\]
defined by
\begin{equation}\label{eq:nK:Dirac}
D (X,f,g)=\left( \tfrac{1}{2}\alpha(dJX)+3X+df+Jdg,d^\ast X+6f, d^\ast(JX)-6g\right). 
\end{equation}

\begin{prop}
$D$ is a self-adjoint elliptic operator.
\proof
Proposition \ref{prop:Dirac:Riemannian} shows that, after choosing appropriate isomorphisms between $\slashed{S}(M)$ and $\Lambda^0\oplus\Lambda^0\oplus\Lambda^1$ we can identify $D$ with $\slashed{D}$ up to a self-adjoint zeroth order term.
\endproof
\end{prop}

\begin{remark}
For every $s\in\R$ one can define a Dirac operator $D_s$ associated with the connection $\nabla ^s = \nabla + \tfrac{s}{2}T$, where $\nabla$ is the Levi-Civita connection and $T=-\Imag\Omega$ is the intrinsic torsion (when $s=1$ this is the so-called \emph{canonical Hermitian connection}). None of these Dirac operators coincides with the operator $D$ defined in \eqref{eq:nK:Dirac}.
\end{remark}

Our interest in the operator $D$ arises from the following identities. Consider the operator
\begin{equation}\label{eq:Dirac+}
D^+\co \Omega^2_{1\oplus 6}\oplus\Omega^4_1 \longrightarrow \Omega^3_{1\oplus 1\oplus 6}
\end{equation}
defined by
\[
(f\omega+X\lrcorner\Real\Omega, \tfrac{1}{2}g\omega^2)\longmapsto \pi_{1\oplus 1\oplus 6}\left( d(f\omega-X\lrcorner\Real\Omega)+\tfrac{1}{2}d^\ast(g\omega^2)\right).
\]
Since by Proposition \ref{prop:differential:2forms} the image of $(f\omega+X\lrcorner\Real\Omega,\tfrac{1}{2}g\omega^2)$ is
\[
\left( \tfrac{1}{2}\alpha(dJX)+df+Jdg + 3X\right)\wedge\omega + \tfrac{1}{2}\left( d^\ast X + 6f\right)\Real\Omega + \tfrac{1}{2}\left( d^\ast JX - 6g\right)\Imag\Omega,
\]  
$D^+$ coincides with $D$ after choosing appropriate identifications of $\Omega^2_{1\oplus 6}\oplus\Omega^4_1$ and $\Omega^3_{1\oplus 1\oplus 6}$ with $\Omega^0\oplus\Omega^0\oplus\Omega^1$. Similarly, the operator
\begin{equation}\label{eq:Dirac-}
D^-\co \Omega^3_{1\oplus1\oplus 6} \longrightarrow \Omega^4_{1\oplus 6}\oplus \Omega^2_1
\end{equation}
defined by
\[
\sigma=JX\wedge\omega-g\Real\Omega+f\Imag\Omega\longmapsto \left( \pi_{1\oplus 6}d\sigma, \pi_1 d^\ast\sigma\right)
\]
coincides with $D$ after identifying $\Omega^4_{1\oplus 6}\oplus \Omega^2_1$ with $\Omega^0\oplus\Omega^0\oplus\Omega^1$ by
\[
(f,g,X)\longmapsto \left( (X\lrcorner\Real\Omega)\wedge\omega - \tfrac{1}{3}f\omega^2, \tfrac{2}{3}g\omega\right).
\]

Much of what follows relies on the mapping properties of the operator $D$ (and therefore of $D^\pm$).

\begin{prop}\label{prop:nK:Dirac} Let $(M,g,\omega,\Omega)$ be a complete nearly K\"ahler $6$--manifold not isometric to the round $6$--sphere. Then the kernel (and cokernel, since $D$ is self-adjont) of $D$ consists of Killing fields that preserve the whole $SU(3)$ structure.
\proof
Suppose that $(X,f,g)$ lies in the kernel of $D$. Then
\begin{itemize}
\item[(a)] $\tfrac{1}{2}\alpha(dJX)+df+Jdg + 3X=0$;
\item[(b)] $d^\ast X+6f=0$;
\item[(c)] $d^\ast (JX)-6g=0$.
\end{itemize}

We apply $d^\ast\circ J$ to the identity (a): using (ii) and (v) in Proposition \ref{prop:differential:2forms} and (c) we obtain $d^\ast dg+18 g=0$ (recall Remark \ref{rmk:T:Tast}!) and therefore $g=0$. Hence
\begin{itemize}
\item[(a')] $\tfrac{1}{2}\alpha(dJX)+df+ 3X=0$;
\item[(b')] $d^\ast X+6f=0$;
\item[(c')] $d^\ast (JX)=0$.
\end{itemize}

Now, (a') and (b') together with \eqref{eq:Dirac:prelim:identities} imply
\[
dJX=-\tfrac{1}{3}(d^\ast X)\,\omega + \tfrac{1}{2}\alpha (dJX)\lrcorner\Real{\Omega} + \pi_8(dJX)=2f\omega-(df+3X)\lrcorner\Real{\Omega}+ \pi_8 (dJX).
\]
Using Proposition \ref{prop:differential:2forms} we differentiate this identity and take the inner product with $\Real\Omega$:
\[
0=\tfrac{1}{4}\langle d^2 (JX),\Real{\Omega}\rangle = 6f+\tfrac{1}{2}d^\ast\left( df+3X \right) = \tfrac{1}{2}(d^\ast df-6f).
\]
By Obata's Theorem $f=0$ under the assumptions of the Proposition.

It remains to show that a vector field $X$ such that $dJX=-3X\lrcorner\Real\Omega+\pi_8(dJX)$ and $d^\ast (JX)=0=d^\ast X$ must preserve the $SU(3)$ structure. Since $d^\ast X=0$ and thus $dJX\wedge\omega^2=0$, by Lemma \ref{lem:identities:2forms}.(i) we have
\[
18\| X\|^2_{L^2}-\| \pi_8(dJX)\| _{L^2}^2=\int_{M}{dJX\wedge dJX\wedge\omega}=3\int_{M}{JX\wedge dJX\wedge\Real\Omega}=18\| X\|^2_{L^2}.
\]
Thus $0=dJX+3X\lrcorner \Real\Omega=\mathcal{L}_X\omega$. Moreover,  by Proposition \ref{prop:differential:2forms}.(iv) $d(X\lrcorner \Real\Omega)\in\Omega^2_{12}$ and therefore the formula of Lemma \ref{lem:Lie:derivative} implies that $\mathcal{L}_X\Real\Omega=d(X\lrcorner \Real\Omega)=0$.
\endproof
\end{prop}

\begin{remark}
When $M$ is the round $6$--sphere with its standard nearly K\"ahler structure the kernel of $D$ consists of elements of the form $(X-\nabla f,f,0)$ where $X$ is a Killing field such that $\mathcal{L}_X\omega=0=\mathcal{L}_X\Real\Omega$ and $f$ satisfies $d^\ast df=6f$.
\end{remark}

\subsection{Hodge decomposition on nearly K\"ahler manifolds}

As immediate corollaries of Proposition \ref{prop:nK:Dirac} we obtain useful decompositions of the space of $3$--forms and $4$--forms. First we fix some notation.

\begin{definition}\label{def:Killing}
We denote by $\mathcal{K}$ the space of Killing fields of $M$. 
\end{definition}
By \cite[Corollary 3.2]{Moroianu:Semmelmann:Unit:Killing} if $M$ is not isometric to the round $6$--sphere then every $X\in\mathcal{K}$ preserves the whole $SU(3)$ structure. Hence by Proposition \ref{prop:nK:Dirac} when $M$ is not isometric to the round $6$--sphere the kernel and cokernel of $D$, and therefore of $D^{\pm}$ in \eqref{eq:Dirac+} and \eqref{eq:Dirac-}, are isomorphic to $\mathcal{K}$.

\begin{prop}\label{prop:Decomposition:3:4:forms} Let $(M,\omega,\Omega)$ be a nearly K\"ahler manifold not isometric to the round $6$--sphere. Then the following holds.
\begin{enumerate}
\item $\Omega^3 = \{ X\wedge\omega\, | \, X \in \mathcal{K}\} \oplus d\Omega^2_{1\oplus6}\oplus d^\ast\Omega^4_1\oplus \Omega^3_{12}$.
\item We have an $L^2$--orthogonal decomposition $\Omega^3_{\text{exact}} = d\Omega^2_{1\oplus6}\oplus \Omega^3_{12,\text{exact}}$.
\item $\Omega^4 = \{ (X\lrcorner\Real\Omega)\wedge\omega \, | \, X\in\mathcal{K}\} \oplus d\Omega^3_{1\oplus 1\oplus6} \oplus \Omega^4_8$. More precisely, for every $4$--form $\sigma$ there exists unique $X\in\mathcal{K}$, $Y\in\mathcal{K}^{\perp_{L^2}}\subset\Omega^1$, $f\in \Omega^0$ and $\sigma_0 \in \Omega^4_{8}$ such that
\[
\sigma = (X\lrcorner\Real\Omega)\wedge\omega + d(JY\wedge\omega + f\Imag\Omega) + \sigma_0.
\]
\item We have an $L^2$--orthogonal decomposition $\Omega^4_{\text{exact}} = d\Omega^3_{1\oplus 1\oplus6} \oplus \Omega^4_{8,\text{exact}}$.
\end{enumerate}
\proof
The identification of $D^+$ with $D$ and Proposition \ref{prop:nK:Dirac} immediately imply (i). In order to deduce (ii) from (i) observe that $2d^\ast (X\wedge\omega)=-\ast\mathcal{L}_X\omega^2=0$ for all $X\in\mathcal{K}$. In particular, $\{ X\wedge\omega \, | \, X\in\mathcal{K}\} \oplus d^\ast\Omega^4_1$ is $L^2$--orthogonal to exact forms (and point-wise orthogonal to $\Omega^3_{12}$). As for the $L^2$--orthogonality statement, observe that by Remark \ref{rmk:Lie:derivatives} if $\rho_0\in\Omega^3_{12}$ is closed then $d^\ast\rho_0\in\Omega^2_8$.

Similarly, the decomposition (iii) follows from the identification of $D^-$ with $D$ and Proposition \ref{prop:nK:Dirac}. Thus every $4$--form $\sigma$ can be written as $\sigma = (X\lrcorner\Real\Omega)\wedge\omega + d(JY\wedge\omega - g\Real\Omega + f \Imag\Omega) + \sigma_0$ with $X\in\mathcal{K}$, $Y \in \mathcal{K}^{\perp_{L^2}}$ and $\sigma_0\in\Omega^4_8$. These constraints uniquely specify $(X,Y,f,g,\sigma_0)$ up to prescribing $d^\ast (JY\wedge\omega - g\Real\Omega + f \Imag\Omega)\wedge\omega^2$. Now, going back to the first step in the proof of Proposition \ref{prop:nK:Dirac} one can show that for every $(Y',f')$ with $Y'\in\mathcal{K}^{\perp_{L^2}}$ there exists a unique choice of $g'$, namely $g'=\frac{1}{3}d^\ast(JY')$, such that every solution $(Y,f,g)$ to $D(Y,f,g)=(Y',f',g')$ satisfies $g=0$.

Finally, (iii) and the fact that $d^\ast\left( (X\lrcorner\Real\Omega)\wedge\omega\right)=-\ast \mathcal{L}_X\Real\Omega=0$ for every $X\in\mathcal{K}$ imply (iv). As for the $L^2$--orthogonality statement, observe that by Proposition \ref{prop:differential:2forms}.(iii) if $\sigma_0\in\Omega^4_{8}$ is closed then $d^\ast\sigma_0\in\Omega^3_{12}$.
\endproof
\end{prop}

Proposition \ref{prop:Decomposition:3:4:forms} plays a crucial role in our treatment of the deformation theory of nearly K\"ahler manifolds. As an additional interesting application, we recover results of Verbitsky \cite[Theorem 6.2 and Remark 6.4]{Verbitsky} on the cohomology of nearly K\"ahler $6$--manifolds. Note that the full information about the cohomology of nearly K\"ahler manifolds is contained in degree $2$ and $3$ since $\pi_1 (M)$ is finite. 

\begin{theorem}\label{thm:Hodge:theory}
Let $(M,\omega,\Omega)$ be a complete nearly K\"ahler manifold. Then every harmonic $2$--form on $M$ lies in $\Omega^2_8$ and every harmonic $3$--form lies in $\Omega^3_{12}$.
\proof
If $M$ is diffeomorphic to $S^6$ the result is vacuous. We can therefore assume that $(M,g)$ is not isometric to the round $6$--sphere.

Let $\eta$ be a closed and coclosed $2$--form. By Proposition \ref{prop:Decomposition:3:4:forms}.(iii) we can write
\[
\eta=X\lrcorner\Real\Omega + d^\ast \rho +\eta_0,
\]
with $X\in\mathcal{K}$, $\rho\in\Omega^3_{1\oplus 1\oplus 6}$ and $\eta_0\in\Omega^2_8$.

Now, on one side $d^\ast\eta=0$ is equivalent to $d^\ast\eta_0=-d^\ast(X\lrcorner\Real\Omega)=6JX$. However,
\[
\langle d^\ast\eta_0, JX\rangle _{L^2}=\langle \eta_0, dJX\rangle _{L^2}=-3\langle \eta_0, X\lrcorner\Real\Omega\rangle_{L^2}=0
\]
for every Killing field $X$. Thus $X=0=d^\ast\eta_0$ and therefore $d\eta_0\in\Omega^3_{12}$ by Proposition \ref{prop:differential:2forms}.(iii). Then
\[
0=\langle \eta, d^\ast\rho \rangle_{L^2}=\| d^\ast\rho\| _{L^2}^2 +\langle \rho,d\eta_0\rangle_{L^2}=\| d^\ast\rho\| _{L^2}^2.
\]

Similarly, let $\rho$ be a closed and coclosed $3$--form and write
\[
\rho = X\wedge\omega + d\eta + d^\ast (f\omega^2)+\rho_0
\]
for $X\in\mathcal{K}$, $\eta\in\Omega^2_{1\oplus 6}$, $f\in C^\infty$ and $\rho_0\in\Omega^3_{12}$. 

We already observed that $X\wedge\omega$ is $L^2$--orthogonal to exact forms whenever $X$ is a Killing field. Furthermore, since in this case $d^\ast JX=0$ we also have
\[
\langle X\wedge\omega, d^\ast(f\omega^2)\rangle _{L^2}=\langle d(X\wedge\omega), f\omega^2\rangle _{L^2}=\langle dX\wedge\omega-3X\wedge\Real\Omega, f\omega^2\rangle _{L^2}=0.
\]
Thus integration by parts using the fact that $\rho$ is closed and coclosed yields
\begin{equation}\label{eq:Hodge:theory:1}
0=\langle \rho, d^\ast (f\omega^2)\rangle _{L^2}=\| d^\ast (f\omega^2)\|^2_{L^2} + \langle d\rho_0, f\omega^2\rangle_{L^2}, \qquad 0=\langle \rho, d\eta\rangle _{L^2}=\| d\eta\|^2_{L^2} + \langle d^\ast\rho_0, \eta\rangle_{L^2}.
\end{equation}

Since $\rho_0\wedge\omega=0$, differentiation yields $d\rho_0\wedge\omega=3\rho_0\wedge\Real\Omega=0$ and therefore the first identity in \eqref{eq:Hodge:theory:1} implies
\[
\rho = X\wedge\omega + d\eta + \rho_0.
\]

In particular $d\rho_0=-d(X\wedge\omega)$. Consider the $\Omega^4_6$--component of this identity. We have
\[
-d(X\wedge\omega)=-dX\wedge\omega+3X\wedge\Real\Omega = -dX\wedge\omega + 3(JX\lrcorner\Real\Omega)\wedge\omega
\]
by \eqref{eq:identities:1:forms:a}. Moreover, by \eqref{eq:dX:dJX} and the fact that $\frac{1}{2}\alpha(dJX)=-3X$ since $X\in\mathcal{K}$ we have
\[
-\pi_6 (dX\wedge\omega) = -\left( \tfrac{1}{2}\alpha(dX)\lrcorner\Real\Omega \right)\wedge\omega = (JX\lrcorner\Real\Omega)\wedge\omega.
\]
Thus $\pi_6 (d\rho_0)=4(JX\lrcorner\Real\Omega)\wedge\omega$. However, by Lemma \ref{lem:Lie:derivative}
\[
\langle d\rho_0, (JX\lrcorner\Real\Omega)\wedge\omega\rangle _{L^2}=\langle d\rho_0, X\wedge\Real\Omega\rangle _{L^2}=\langle \rho_0, d(X\lrcorner\Real\Omega)\rangle _{L^2}=0
\]
since $X\in\mathcal{K}$. We conclude that $X=0$ and $\rho=d\eta+\rho_0$.

In particular, $d\rho_0=0$. Now, on one side $\ast\left( d^\ast\rho_0\wedge\omega^2\right)=-2\ast\left( d\ast\rho_0\wedge\omega\right)=0$ since $\ast\rho_0\in\Omega^3_{12}$. On the other side, $d\rho_0=0$ and \eqref{eq:alpha:primitive:(2,1)} imply that $d\ast\rho_0$ has no $\Omega^4_6$ component. We conclude that $d^\ast\rho_0\in\Omega^2_8$ and therefore $\langle d^\ast\rho_0, \eta\rangle_{L^2}=0$. In view of the second identity in \eqref{eq:Hodge:theory:1} the Theorem is proved. 
\endproof
\end{theorem}

\begin{remark}
Theorem \ref{thm:Hodge:theory} has an interesting application, suggested to us by Gon\c calo Oliveira, to gauge theory on nearly K\"ahler manifolds. On every nearly K\"ahler manifold one can define the notion of a \emph{pseudo-Hermitian-Yang-Mills} connection: a connection $A$ on a principal $G$--bundle whose curvature $F_A$ is a primitive $(1,1)$--form with values in the adjoint bundle. The interest in pseudo-Hermitian-Yang-Mills connections on a nearly K\"ahler manifold $M$ arises from the fact that they correspond to scale-invariant $G_2$--instantons on the Riemannian cone over $M$. In particular, pseudo-Hermitian-Yang-Mills connections on the round $S^6$ model isolated singularities of $G_2$ instantons on smooth $G_2$ manifolds. Now, Theorem \ref{thm:Hodge:theory} implies that every line bundle $L$ over a nearly K\"ahler manifold admits a pseudo-Hermitian-Yang-Mills connection. Indeed, the connection with curvature given by the harmonic representative of $-2\pi i\, c_1 (L)$ is pseudo-Hermitian-Yang-Mills.
\end{remark}

\section{Deformations of nearly K\"ahler manifolds}

Let $(M,\omega,\Omega)$ be a nearly K\"ahler manifold not isometric to the round $6$--sphere. We are going to study the problem of deforming $(\omega,\Omega)$ to a nearby nearly K\"ahler structure $(\omega',\Omega')$.

This will be done in two steps. First we use Proposition \ref{prop:Decomposition:3:4:forms} to define a slice for the action of the diffeomorphism group on the space of $SU(3)$ structures close to $(\omega,\Omega)$. This choice of slice will allow us to give a streamlined proof of the identification of the space of infinitesimal nearly K\"ahler deformations with an eigenspace of the Laplacian acting on coclosed primitive $(1,1)$ forms, a result due to Moroianu--Nagy--Semmelmann \cite{Moroianu:Nagy:Semmelmann}.

In order to proceed beyond the linearised level, we find convenient to exploit some observations due to Hitchin \cite{Hitchin} to enlarge the spaces under considerations: we will reinterpret the nearly K\"ahler equations \eqref{eq:NK} as the vanishing of a smooth map on the space of pairs of an exact stable $4$--form and an exact stable $3$--form, without requiring a priori any compatibility condition. Studying the mapping properties of the linearisation of this map will identify possible obstructions to integrate infinitesimal nearly K\"ahler deformations to genuine nearly K\"ahler structures.   

\subsection{Deformations of $SU(3)$ structures}

Let $M$ be a $6$--manifold. We denote by $\mathcal{C}$ the space of $SU(3)$ structures on $M$, \ie the set of all $(\Real\Omega,\Imag\Omega,\tfrac{1}{2}\omega^2,\omega)\in\Omega^3(M)\times\Omega^3(M)\times\Omega^4(M)\times\Omega^2(M)$ such that $\Real\Omega$ is a stable $3$--form with dual $\Imag\Omega$, $\tfrac{1}{2}\omega^2$ is a stable $4$--form with dual $\omega$ and the constraints $\omega\wedge\Real\Omega=0=4\omega^3-3\Real\Omega\wedge\Imag\Omega$ are satisfied. We will label an $SU(3)$ structure by $(\Real\Omega,\tfrac{1}{2}\omega^2)$.

In fact, for the purposes of doing analysis it will be more appropriate to introduce H\"older spaces and consider the subspace $\mathcal{C}^{k,\alpha}$ of $SU(3)$ structures such that $\omega,\Omega$ are of class $C^{k,\alpha}$ for some $k\geq 1$ and $\alpha\in (0,1)$. Thus $\mathcal{C}^{k,\alpha}$ is continuously embedded as a submanifold of the space of differential forms of class $C^{k,\alpha}$. For ease of notation we will drop the index $^{k,\alpha}$ when this will not be essential.

Given a point $\Lie{c}=(\Real\Omega,\tfrac{1}{2}\omega^2)\in\mathcal{C}$ we now define the tangent space $T_{\Lie{c}}\mathcal{C}$ and an exponential map $\text{exp}\co \mathcal{U}\ra\mathcal{C}$ from a sufficiently small neighbourhood $\mathcal{U}$ of the origin in $T_{\Lie{c}}\mathcal{C}$.

\begin{lemma}\label{lem:Tangent:Space}
$T_{\Lie{c}}\mathcal{C}$ is the set of all $(\rho,\hat{\rho},\sigma,\hat{\sigma})\in \Omega^3(M)\times\Omega^3(M)\times\Omega^4(M)\times\Omega^2(M)$ such that
\[
\begin{gathered}
\rho = -X\wedge\omega + 3\lambda\Real\Omega + 3\mu\Imag\Omega +\rho_0,\\
\hat{\rho}= -JX\wedge\omega -3\mu\Real\Omega + 3\lambda \Imag\Omega -\ast\rho_0,\\
\sigma = \hat{\sigma}\wedge\omega, \qquad \hat{\sigma}=X\lrcorner\Real\Omega + 2\lambda\omega +\eta_0,
\end{gathered}
\]
for some $X\in\mathcal{X}(M)$, $\lambda,\mu\in\Omega^0(M)$, $\rho_0\in\Omega^3_{12}$, $\eta_0\in\Omega^2_8$.
\proof
The Lemma is a straightforward consequence of the linearisation of the constraints \eqref{eq:SU(3):structure:Constraints} and Proposition \ref{prop:Linearisation:Hitchin:dual}.
\endproof
\end{lemma}

Note that linear $SU(3)$ structures on $\R^6$ are parametrised by $GL(6,\R)/SU(3)$. Following \cite{Goto}, given a manifold $M$ endowed with an $SU(3)$ structure $\Lie{c}=(\Real\Omega,\tfrac{1}{2}\omega^2)$ we define the sub-bundle $SU(TM)$ of the frame bundle $GL(TM)$ of $M$ whose sections are the bundle isomorphisms of $TM$ which preserve $(\Real\Omega,\tfrac{1}{2}\omega^2)$ at each point. Then $\mathcal{C}\simeq GL(TM)/SU(TM)$ under the isomorphism that associate to each $g\in GL(TM)$ the $SU(3)$ structure $g_\ast(\Real\Omega,\tfrac{1}{2}\omega^2)$.

From this perspective $T_{\Lie{c}}\mathcal{C}$ is given by elements $A_\ast(\Real\Omega,\tfrac{1}{2}\omega^2)$ for $A\in GL(M)\times_{GL(6,\R)}\Lie{m}$. Here $\Lie{m}$ is a complement of $\Lie{su}(3)$ in $\Lie{gl}(6,\R)$ and $GL(M)\times_{GL(6,\R)}\Lie{m}$ is the bundle associated with the representation $\Lie{m}$ of $GL(6,\R)$. For a sufficiently small neighbourhood $\mathcal{U}$ of the origin in $T_{\Lie{c}}\mathcal{C}$, these identifications and the exponential map $\Lie{m}\ra GL(6,\R)$ induce an exponential map
\begin{equation}\label{eq:Exponential}
\text{exp}\co \mathcal{U} \longrightarrow \mathcal{C}, \qquad \text{exp}(A_\ast\Lie{c})=\sum_{n=0}^\infty{\tfrac{1}{n!}(A_\ast)^n\Lie{c}}=\Lie{c} + A_\ast\Lie{c}+\tfrac{1}{2}A_\ast A_\ast\Lie{c}+\dots
\end{equation}
which is a homeomorphism onto its image.
 
\subsection{Slice to the action of the diffeomorphism group}

The first step in studying the deformation theory of nearly K\"ahler structures is to find a slice for the action of the diffeomorphism group.

Let $\Lie{c}=(\Real\Omega,\tfrac{1}{2}\omega^2)$ be a nearly K\"ahler structure such that the induced metric is not isometric to the round metric on $S^6$. Denote by $\mathcal{O}_{\Lie{c}}$ the orbit of $\Lie{c}$ under the action of the group $\text{Diff}_0^{k+1,\alpha}(M)$ of $C^{k+1,\alpha}$--diffeomorphisms of $M$ isotopic to the identity. The tangent space $T_{\Lie{c}}\mathcal{O}_{\Lie{c}}$ is the space of Lie derivatives $\mathcal{L}_X (\Real\Omega,\tfrac{1}{2}\omega^2)$ for $X\in C^{k,\alpha}(TM)$.

We now use Proposition \ref{prop:Decomposition:3:4:forms} to define a complement $\mathcal{W}$ of $T_{\Lie{c}}\mathcal{O}_{\Lie{c}}$ in $T_{\Lie{c}}\mathcal{C}$. Given $(\rho,\hat{\rho},\sigma,\hat{\sigma})\in T_{\Lie{c}}\mathcal{C}$, write
\[
\sigma = -(X\lrcorner\Real\Omega)\wedge\omega + d(JY\wedge\omega - f\Imag\Omega) + \sigma_0
\]
for unique $X\in\mathcal{K}$, $Y\in\mathcal{K}^{\perp_{L^2}}$, $f\in \Omega^0$ and $\sigma_0 \in \Omega^2_{8}$. The term $d(JY\wedge\omega)=\mathcal{L}_Y\left( \tfrac{1}{2}\omega^2\right)$ is a Lie derivative. Thus up to an element in $T_{\Lie{c}}\mathcal{O}_{\Lie{c}}$ we can assume that
\begin{subequations}\label{eq:Slice}
\begin{equation}
\sigma = (-X+\nabla f+J\nabla g)\lrcorner\Real\Omega\wedge\omega + 2f\omega^2 + \sigma_0.
\end{equation}
Then Lemma \ref{lem:Tangent:Space} forces
\begin{equation}
\rho = (X+df-Jdg)\wedge\omega + 3f\Real\Omega+3g\Imag\Omega + \rho_0.
\end{equation}
\end{subequations}
for some $\rho_0\in\Omega^3_{12}$. $\hat{\rho}$, $\hat{\sigma}$ are then determined by $\rho$ and $\sigma$ via Proposition \ref{prop:Linearisation:Hitchin:dual}.

Hence Proposition \ref{prop:Decomposition:3:4:forms} yields a splitting $T_{\Lie{c}}\mathcal{C}=T_{\Lie{c}}\mathcal{O}_{\Lie{c}}\oplus \mathcal{W}$ where $\mathcal{W}\simeq \mathcal{K}\times\Omega^0(M)\times\Omega^0(M)\times\Omega^2_{8}(M)\times\Omega^3_{12}(M)$ is the space of all $(\rho,\hat{\rho},\sigma,\hat{\sigma}) \in T_{\Lie{c}}\mathcal{C}$ with $\rho$ and $\sigma$ of the form \eqref{eq:Slice}. 

\begin{prop}[{\cite[Theorems 3.1.4 and 3.1.7]{Nordstrom:Thesis}}]\label{prop:Slice}
There exists an open neighbourhood $\mathcal{U}\subset \mathcal{W}$ of the origin such that $\mathcal{S}=\text{exp}(\mathcal{U})$ is a slice for the action of $\text{Diff}_0^{k+1,\alpha}(M)$ on a sufficiently small neighbourhood of $\Lie{c}=(\Real\Omega,\tfrac{1}{2}\omega^2)\in\mathcal{C}$.
\end{prop}

\subsection{Infinitesimal deformations}

With this choice of slice we can easily determine infinitesimal deformations of the nearly K\"ahler structure $\Lie{c}$, thus recovering the result of Morianu--Nagy--Semmelmann \cite[Theorem 4.2]{Moroianu:Nagy:Semmelmann} as a simple consequence of Proposition \ref{prop:Decomposition:3:4:forms}.

\begin{theorem}\label{thm:Infinitesimal:deformations}
Let $(M,\omega,\Omega)$ be a nearly K\"ahler manifold non-isometric to the round $6$--sphere. Then infinitesimal deformations of the nearly K\"ahler structure modulo diffeomorphisms are in one to one correspondence with pairs $(\rho_0,\sigma_0)\in\Omega^3_{12,\text{exact}}\oplus\Omega^4_{8,\text{exact}}$ satisfying
\begin{equation}\label{eq:NK:Linearisation:Gauged}
-d\ast\sigma_0-3\rho_0=0, \qquad -d\ast\rho_0+4\sigma_0=0.
\end{equation}
\proof
Linearising the nearly K\"ahler equations \eqref{eq:NK} we see that an infinitesimal deformations $(\rho,\hat{\rho},\sigma,\hat{\sigma})\in T_{\Lie{c}}\mathcal{C}$ of $\Lie{c}$ as an $SU(3)$ structure is an infinitesimal nearly K\"ahler deformation if
\begin{equation}\label{eq:NK:Linearisation}
d\hat{\sigma}-3\rho=0=d\hat{\rho}+4\sigma.
\end{equation}
In particular, $(\rho,\sigma)\in \Omega^3_{\text{exact}}\oplus\Omega^4_{\text{exact}}$. Thus if we further assume that $(\rho,\hat{\rho},\sigma,\hat{\sigma})\in \mathcal{W}$ then by Proposition \ref{prop:Decomposition:3:4:forms}.(iv) we have
\[
\sigma = -d\left( f\Imag\Omega\right) + \sigma_0
\]
with $\sigma_0\in\Omega^4_{8,\text{exact}}$. 
Lemma \ref{lem:Tangent:Space} then implies that
\[
\rho = d(f\omega) + \rho_0
\]
for some $\rho_0\in\Omega^3_{12,\text{exact}}$ and $\hat{\sigma} = 2f\omega-\nabla f\lrcorner\Real\Omega -\ast \sigma_0$.

Consider the equation $d\hat{\sigma}-3\rho=0$: by Proposition \ref{prop:differential:2forms}.(ii) and (iv)
\[
0=\tfrac{1}{4}\langle d\hat{\sigma}-3\rho,\Real\Omega\rangle = \tfrac{1}{2}\left( d^\ast df-6f\right)
\]
and therefore $f=0$ by Obata's Theorem.
\endproof
\end{theorem}

As noted in \cite[Theorem 4.2]{Moroianu:Nagy:Semmelmann}, a solution $(\rho_0,\sigma_0)$ of \eqref{eq:NK:Linearisation:Gauged} is uniquely determined by the coclosed primitive $(1,1)$--form $\ast\sigma_0 \in \Omega^2_8$. This satisfies $\triangle (\ast \sigma_0)=12(\ast\sigma_0)$. The space of infinitesimal nearly K\"ahler deformations is therefore identified with the intersection of the space of coclosed primitive $(1,1)$--forms with the eigenspace of eigenvalue $12$ of the Laplacian acting on $2$--forms.

\subsection{The deformation problem}

To proceed further we would like to reformulate the nearly K\"ahler equations \eqref{eq:NK} as the vanishing of a smooth map $\Phi$ on the slice $\mathcal{S}$. In fact, it will be convenient to enlarge the space $\mathcal{C}$ of $SU(3)$ structures dropping the requirement that the constraints \eqref{eq:SU(3):structure:Constraints} are satisfied. This approach to the deformation theory of nearly K\"ahler manifolds is very natural from the point of view introduced in Hitchin's seminal paper \cite{Hitchin}.

Hitchin's first observation is that if $\Real\Omega$ is a stable $3$--form with dual $\Imag\Omega$ and $\tfrac{1}{2}\omega^2$ is a stable $4$--form with dual $\omega$ satisfying the nearly K\"ahler equations \eqref{eq:NK} then $(\Real\Omega,\tfrac{1}{2}\omega^2)$ defines an $SU(3)$ structure, \ie the compatibility constraints \eqref{eq:SU(3):structure:Constraints} are automatically satisfied. Indeed, observe that
\[
3\,\omega\wedge\Real\Omega = \omega\wedge d\omega=d\left( \tfrac{1}{2}\omega^2\right)=0.
\]
Since $\Imag\Omega$ is the dual of $\Real\Omega$ then $\omega\wedge\Imag\Omega=0$ also. Differentiating this identity we obtain
\[
0=d\omega\wedge\Imag\Omega + \omega\wedge d\Imag\Omega = 3\Real\Omega\wedge\Imag\Omega - 2\omega^3.
\]
Avoiding to impose the constraints \eqref{eq:SU(3):structure:Constraints} from the start let us gain a vector field and a function as additional free parameters.

The appearance of a second vector field as an additional parameter follows from the action of the diffeomorphism group. In order to explain this last point, following Hitchin \cite{Hitchin} we interpret nearly K\"ahler structures as critical points of a diffeomorphism-invariant functional on an open set in the infinite dimensional symplectic vector space $\Omega^3_{\text{exact}}\oplus\Omega^4_{\text{exact}}$.

Let $\widetilde{\mathcal{C}}$ be the space of stable forms $(\Real\Omega,\tfrac{1}{2}\omega^2)\in\Omega^3_{\text{exact}}\times\Omega^4_{\text{exact}}$. We denote by $\Imag\Omega$ and $\omega$ the respective duals. For each element in $\widetilde{\mathcal{C}}$ Hitchin defines volume functionals $V_3(\Real\Omega)$ and $V_4(\tfrac{1}{2}\omega^2)$ such that for any $(\rho,\sigma)\in \Omega^3_{\text{exact}}\oplus\Omega^4_{\text{exact}}$
\[
\delta V_3 (\rho)=-\int{\Imag\Omega\wedge\rho}, \qquad \delta V_4 (\sigma)=\int{\omega\wedge\sigma}.
\]
Furthermore, consider the natural pairing $P$ on $\Omega^3_{\text{exact}}\times\Omega^4_{\text{exact}}$ defined by
\[
P(\rho,\sigma)=\int{\alpha\wedge\sigma}=-\int{\rho\wedge\beta}
\]
for $d\alpha=\rho$ and $d\beta=\sigma$. Now define a functional $\mathcal{E}\co \widetilde{\mathcal{C}}\ra \R$ by
\[
\mathcal{E}(\Real\Omega,\tfrac{1}{2}\omega^2)=3V_3(\Real\Omega)+4V_4(\tfrac{1}{2}\omega^2)-12P(\Real\Omega,\tfrac{1}{2}\omega^2).
\]
Then for every $(\rho,\sigma)\in \Omega^3_{\text{exact}}\times\Omega^4_{\text{exact}}$ we have
\[
\delta\mathcal{E}(\rho,\sigma)=-3\int{\left( \Imag\Omega+4\beta\right)\wedge\rho} + 4 \int{\left( \omega-3\alpha\right)\wedge\sigma}
\]
where $d\alpha=\Real\Omega$ and $d\beta=\tfrac{1}{2}\omega^2$. Critical points of $\mathcal{E}$ are therefore pairs $(\Real\Omega,\tfrac{1}{2}\omega^2)\in\widetilde{\mathcal{C}}$ satisfying $d\Imag\Omega+2\omega^2=0=d\omega-3\Real\Omega=0$, \ie nearly K\"ahler structures.

Since $\mathcal{E}$ is diffeomorphism invariant it is clear that $\delta\mathcal{E}$ vanishes in the direction of Lie derivatives at any point in $\widetilde{\mathcal{C}}$ and this explain the freedom to throw in the equations an extra vector field.

In the following Proposition we take stock of these observations to rewrite the nearly K\"ahler equations as a the vanishing of a smooth map.

\begin{prop}\label{prop:NK:revisited}
Suppose that $(\Real\Omega,\tfrac{1}{2}\omega^2)\in\widetilde{\mathcal{C}}$ satisfies
\[
d\omega-3\Real\Omega=0, \qquad d\Imag\Omega+2\omega^2=d\ast d(Z\lrcorner\Real\Omega)
\]
for some vector field $Z$. Here the Hodge--$\ast$ is associated with a fixed background metric.

Then $(\Real\Omega,\tfrac{1}{2}\omega^2)$ defines a nearly K\"ahler structure.
\proof
It is enough to show that $d(Z\lrcorner\Real\Omega)=0$.

First of all observe as above that $d\omega-3\Real\Omega=0=d\omega\wedge\omega$ imply that $\omega\wedge\Real\Omega=0$ and therefore $(Z\lrcorner\Real\Omega)\wedge\omega^2=0$.

Moreover, since $d\Real\Omega=0$ by Proposition \ref{prop:Torsion:SU(3):structures} we also have $\int{(Y\lrcorner\Real\Omega)\wedge d\Imag\Omega}=0$ for every vector field $Y$. Indeed, we can conformally rescale $\omega$ so that both constraints \eqref{eq:SU(3):structure:Constraints} are satisfied. Denote by $\widetilde{\omega}$ this rescaled form. Then $(\widetilde{\omega},\Omega)$ is an $SU(3)$ structure. In particular, Proposition \ref{prop:Torsion:SU(3):structures} can now be applied: $d\Imag\Omega$ has no component of type $\Omega^4_6$ with respect to $(\widetilde{\omega},\Omega)$. Finally, note that $\int{(Y\lrcorner\Real\Omega)\wedge d\Imag\Omega}=\langle (Y\lrcorner\Real\Omega)\wedge\widetilde{\omega}, d\Imag\Omega \rangle_{L^2}$, where the $L^2$ inner product is computed using the metric induced by $(\widetilde{\omega},\Omega)$.

Integrating by parts we therefore obtain
\[
\| d(Z\lrcorner\Real\Omega)\| ^2_{L^2} = \int{ (Z\lrcorner\Real\Omega)\wedge (d\Imag\Omega+2\omega^2) }=0.\qedhere
\]
\end{prop}

Fix a nearly K\"ahler structure $(\Real\Omega,\tfrac{1}{2}\omega^2)$ and assume that the induced metric in not isometric to the round metric on $S^6$. For every $\rho\in\Omega^3_{exact}$ and $\sigma\in\Omega^4_{exact}$ sufficiently small in $C^{k,\alpha}$ the forms $\Real\Omega'=\Real\Omega+\rho$ and $\tfrac{1}{2}\omega^2+\sigma$ are still stable forms. We therefore have a ``linear'' exponential map $\widetilde{\text{exp}}$ from a neighbourhood of the origin in $\Omega^3_{exact}\times \Omega^4_{exact}$ into a sufficiently small neighbourhood of $(\Real\Omega,\tfrac{1}{2}\omega^2)$ in $\widetilde{\mathcal{C}}$.

Proposition \ref{prop:Decomposition:3:4:forms} could be used to define a slice to the action of the diffeomorphism group also in this case (but we won't really need this beyond the tangent space level): consider the image under $\widetilde{\text{exp}}$ of a small neighbourhood $\mathcal{U}$ of the origin in the subspace $\widetilde{\mathcal{W}}$ of $\Omega^3_{exact}\times \Omega^4_{exact}$ consisting of pairs $(\rho, \sigma)$ where $\sigma = d(g\Imag\Omega) + \sigma_0$ with $\sigma_0\in\Omega^4_{8,exact}$.  

Now define a map
\begin{equation}\label{eq:NK:revisited}
\begin{gathered}
\Phi\co \mathcal{U} \times C^{k+1,\alpha}(TM) \longrightarrow C^{k-1,\alpha}(\Lambda^3 T^\ast M\oplus \Lambda^4 T^\ast M ),\\
\Phi (\rho,\sigma,Z)=\left( d\omega' - 3\Real\Omega', d\Imag\Omega' + 2\omega'\wedge\omega' + d\ast d(Z\lrcorner\Real\Omega')\right).
\end{gathered}
\end{equation}
Here $\Imag\Omega'$ is the dual of $\Real\Omega'=\Real\Omega+\rho$, $\omega'$ is the dual of $\tfrac{1}{2}\omega^2 + \sigma$ and the Hodge--$\ast$ is computed with respect to the metric induced by the  nearly K\"ahler structure $(\Real\Omega,\tfrac{1}{2}\omega^2)$. By Proposition \ref{prop:NK:revisited} nearly K\"ahler structures close to $(\Real\Omega,\tfrac{1}{2}\omega^2)$ are parametrised modulo diffeomorphisms by the zero locus of $\Phi$. 

\subsection{Obstructions}

The linearisation $D\Phi$ of \eqref{eq:NK:revisited} at $(\Real\Omega,\tfrac{1}{2}\omega^2,0)$ is defined by
\begin{equation}\label{eq:NK:Linearisation:revisited}
d\hat{\sigma}-3\rho, \qquad d\hat{\rho}+4\sigma + d\ast d(Z\lrcorner\Real\Omega),
\end{equation}
where $(\rho,\sigma)\in \mathcal{\widetilde{W}}$ and $\hat{\rho}, \hat{\sigma}$ are their images under the linearisation of Hitchin's duality map in Proposition \ref{prop:Linearisation:Hitchin:dual}. Using Proposition \ref{prop:Decomposition:3:4:forms} we can write:
\begin{subequations}\label{eq:Tangent:vector}
\begin{equation}
\rho = d(f\omega + X\lrcorner\Real\Omega) + \rho_0
\end{equation}
with $\rho_0\in\Omega^3_{12,exact}$ and $X\in\mathcal{K}^{\perp_{L^2}}$ and
\begin{equation}
\hat{\rho}=\ast d(f\omega ) -4JX\wedge\omega -\ast\rho_0 +d(X\lrcorner\Imag\Omega)=-Jdf\wedge\omega +3f\Imag\Omega -4JX\wedge\omega -\ast\rho_0 +d(X\lrcorner\Imag\Omega), 
\end{equation}
(here we used the fact that $d(X\lrcorner\Real\Omega)=\mathcal{L}_X\Real\Omega$ to calculate the corresponding component $\mathcal{L}_X\Imag\Omega$ of $\hat{\rho}$);
\begin{equation}
\sigma = -d(g\Imag\Omega\omega) + \sigma_0,
\end{equation}
\begin{equation}
\hat{\sigma} = -\ast d(g\Imag\Omega) -2g\omega -\ast \sigma_0= -\nabla g\lrcorner\Real\Omega +2g\omega -\ast \sigma_0.
\end{equation}
\end{subequations}
 
In order to describe the zero locus of $\Phi$ using the Implicit Function Theorem we need to study the mapping properties of $D\Phi$, \ie given $(\alpha,\beta)\in \Omega^3_{exact}(M)\times\Omega^4_{exact}(M)$ we need to study the equation
\[
D\Phi(\rho,\sigma,Z)=(\alpha,\beta).
\]

\begin{prop}\label{prop:Obstructions}
A pair $(\alpha,\beta)\in\Omega^3_{exact}\times\Omega^4_{exact}$ lies in the image of $D\Phi$ if and only if
\[
\langle d^\ast \alpha + 3\ast \beta, \xi\rangle _{L^2} =0
\]
for every co-closed primitive $(1,1)$ form $\xi$ such that $\triangle\xi=12\xi$.
\proof
We have to solve
\begin{equation}\label{eq:NK:Linearised:revisited}
d\hat{\sigma}-3\rho = \alpha, \qquad d\hat{\rho}+4\sigma +d\ast d(Z\lrcorner\Real\Omega) = \beta
\end{equation}
for $\rho,\alpha\in\Omega^3_{exact}$, $\sigma,\beta\in \Omega^4_{exact}$ and a vector field $Z$.

Using \eqref{eq:Tangent:vector} we calculate
\[
\begin{gathered}
d\hat{\sigma}-3\rho=d\Big( -(3X+\nabla g)\lrcorner\Real\Omega +(2g-3f)\omega  \Big) -(d\ast\sigma_0+3\rho_0),\\
d\hat{\rho}+4\sigma = d\Big( (-4JX+J\nabla f)\wedge\omega + (3f-4g)\Imag\Omega \Big) +(-d\ast\rho_0 + 4\sigma_0).
\end{gathered}
\]

Now set $u=2g-3f$, $v=-3f+4g$ and $U=-3X-\nabla g$. Note that the map $(f,g)\mapsto (u,v)$ is invertible and observe that $-4X+\nabla f = \tfrac{4}{3}(U-\nabla u) + \nabla v$. Thus
\[
\begin{gathered}
d\hat{\sigma}-3\rho=d\Big( U\lrcorner\Real\Omega +u\omega  \Big) -(d\ast\sigma_0+3\rho_0),\\
d\hat{\rho}+4\sigma = d\Big( \tfrac{1}{3}J(4U-4\nabla u+3\nabla v)\wedge\omega - v\Imag\Omega \Big) +(-d\ast\rho_0 + 4\sigma_0).
\end{gathered}
\]

Moreover, by Lemma \ref{lem:Lie:derivative}
\[
\ast d(Z\lrcorner\Real\Omega) = d\Big( JZ\lrcorner\Real\Omega + \tfrac{1}{3}(d^\ast JZ)\omega \Big) - J\Big( \alpha(dJZ)+2Z-\tfrac{1}{3}J\nabla d^\ast(JZ)\Big)\wedge\omega -(d^\ast Z)\Imag\Omega.
\]

Use Proposition \ref{prop:Decomposition:3:4:forms} to write $\alpha = d\eta + \alpha_0$ and $\beta = d\zeta+\beta_0$ with $\eta\in\Omega^2_{1\oplus 6}$, $\alpha_0\in\Omega^2_{8,exact}$, $\zeta\in\Omega^3_{1\oplus 1\oplus 6}$ with $\zeta\wedge\Imag\Omega=0$ and $\beta_0\in \Omega^3_{12,exact}$.

We first look for a solution to \eqref{eq:NK:Linearised:revisited} assuming $\alpha_0=0=\beta_0$. 

Observe that $U$ and $u$ are uniquely determined by $\eta$ up to an element of $\mathcal{K}$. Thus we reduced to study the mapping properties of the operator
\[
D'\co (Z,v)\mapsto \Big( \alpha (dJZ)+2Z-\tfrac{1}{3}J\nabla d^\ast(JZ)-\nabla v, d^\ast Z + v\Big)
\]
of $\Omega^1\times\Omega^0$ into itself. By changing $v$ into $-\tfrac{1}{2}v$ and setting $w=\tfrac{1}{6}d^\ast (JZ)$ the equation $D'(Z,v)=(2X_0,\mu_0)$ can be rewritten as
\[
\tfrac{1}{2}\alpha (dJZ)+Z+dv+Jdw=X_0, \qquad d^\ast Z-2v=\mu_0, \qquad d^\ast JZ -6w=0.
\]

By Proposition \ref{prop:Dirac:Riemannian} $D'$ can therefore be identified with a self-adjoint operator which coincides with $\slashed{D}$ up to a zeroth order term. The same arguments as in the proof of Proposition \ref{prop:nK:Dirac} show that $D'$ has trivial kernel (and cokernel). Indeed, suppose that 
\[
\tfrac{1}{2}\alpha (dJZ)+Z+dv+Jdw=0, \qquad d^\ast Z-2v=0, \qquad d^\ast JZ -6w=0.
\]
Applying $-d^\ast\circ J$ to the first equation and using the third we find $0=d^\ast dw -d^\ast (JZ^\sharp) = d^\ast dw + 6w$ and therefore $w=0$. Thus the first two equations now imply
\[
dJZ=-\tfrac{2}{3}v\omega -(Z+\nabla v)\lrcorner\Real\Omega + \pi_8 (dJZ).
\]
In particular, (using $|\Real\Omega|^2=4$)
\[
0=\tfrac{1}{2}\langle d(dJZ),\Real\Omega\rangle = -4 v + d^\ast Z + d^\ast dv=d^\ast dv -2v
\]
and therefore $v=0$ also. Finally, since $d^\ast Z=0$, by Lemma \ref{lem:identities:2forms}.(i) we have
\[
2\| Z\|^2_{L^2}-\| \pi_8(dJX)\| _{L^2}^2=\int_{M}{dJZ\wedge dJZ\wedge\omega}=3\int_{M}{JZ\wedge dJZ\wedge\Real\Omega}=6\| Z\|^2_{L^2}
\]
and therefore $Z=0$.

We have therefore reduced to solve the equation
\[
-d\ast\sigma_0-3\rho_0 = \alpha_0, \qquad -d\ast\rho_0 + 4\sigma_0=\beta_0,
\]
for $\alpha_0\in\Omega^3_{12,exact}$ and $\beta_0\in\Omega^4_{8,exact}$.

Now, let us rewrite this system of equations as a second order PDE for the coclosed primitive $(1,1)$--form $\hat{\sigma}_0=-\ast\sigma_0$: taking $d^\ast$ of the first equation and using $\ast$ of the second one we find
\[
d^\ast d\hat{\sigma}_0-12\hat{\sigma}_0=d^\ast\alpha_0 + 3\ast\beta_0.
\]
Conversely, given a solution $\hat{\sigma}_0$ of this equation, set $3\rho_0=d\hat{\sigma}_0-\alpha_0 \in\Omega^3_{12,exact}$ to get a solution of the first order system.

It is now clear that a solution of \eqref{eq:NK:Linearised:revisited} exists if and only if $d^\ast\alpha_0 + 3\ast\beta_0$ is $L^2$--orthogonal to the space $\mathcal{O}$ of primitive coclosed $(1,1)$--forms which are eigenforms for the Laplacian with eigenvalue $12$. Furthermore, observe that for every $\eta\in\Omega^2_{1\oplus 6}$ and $\zeta\in\Omega^3_{1\oplus 1\oplus 6}$
\[
\langle d^\ast d\eta, \xi\rangle_{L^2}=12 \langle \eta, \xi\rangle _{L^2}=0, \qquad \langle \ast d\zeta, \xi\rangle_{L^2}=- \langle \zeta, \ast d\xi\rangle _{L^2}=0
\]
for every $\xi\in\mathcal{O}$ since $d\xi\in\Omega^3_{12}$. We therefore conclude that \eqref{eq:NK:Linearised:revisited} is solvable if and only if $d^\ast \alpha + 3\ast \beta$ is $L^2$--orthogonal to $\mathcal{O}$. \endproof
\end{prop}

Proposition \ref{prop:Obstructions} strongly suggests that the deformation theory of nearly K\"ahler manifolds is in general obstructed. In the next section we study a concrete example that shows that this is indeed the case. 

\section{Deformations of the flag manifold}

In this final section we study in details an example. In \cite{Moroianu:Semmelmann} Moroianu--Nagy--Semmelmann study infinitesimal deformations of the homogeneous nearly K\"ahler manifolds (with the exclusion of $S^6$ with its standard nearly K\"ahler structure that has to be considered separately, \cf \cite[Theorem 5.1]{Moroianu:Nagy:Semmelmann}). They found that $S^3\times S^3$ and $\CP^3$ are rigid while the flag manifold $\mathbb{F}_3 = SU(3)/T^2$ admits an $8$--dimensional space of infinitesimal nearly K\"ahler deformations. It has been an open problem to decide whether any of these deformations can be integrated to produce new examples of complete nearly K\"ahler manifolds. Here we apply the framework developed in the previous section to show that all the infinitesimal deformations of the flag manifolds are obstructed.

\subsection{Second order deformations}

Inspired by work of Koiso \cite{Koiso}, we find convenient to introduce the notion of second order deformations and show that infinitesimal deformations of the flag manifold are obstructed to second order.

\begin{definition}\label{def:2nd:order:NK:deformation}
Given a nearly K\"ahler structure $(\Real\Omega_0,\tfrac{1}{2}\omega_0^2)$ and an infinitesimal nearly K\"ahler deformation $(\rho_1,\sigma_1)$ a \emph{second order deformation} in the direction of $(\rho_1,\sigma_1)$ is a $4$--tuple $(\rho_2,\rho'_2,\sigma_2,\sigma'_2)\in \Omega^3(M)\times\Omega^3 (M)\times \Omega^4(M)\times \Omega^2(M)$ such that
\[
\begin{gathered}
\Real\Omega = \Real\Omega_0 + \epsilon \rho_1 + \tfrac{1}{2}\rho_2 , \qquad \Imag\Omega = \Imag\Omega_0 + \hat{\rho}_1 + \tfrac{1}{2}\epsilon^2 \rho'_2,\\
\tfrac{1}{2}\omega ^2 = \tfrac{1}{2}\omega_0^2 + \epsilon \sigma_1 + \tfrac{1}{2}\epsilon^2 \sigma_2, \qquad \omega = \omega_0 + \epsilon \hat{\sigma}_1 + \tfrac{1}{2}\epsilon^2\sigma'_2,
\end{gathered}
\]
is a nearly K\"ahler structure up to terms of order $o(\epsilon^2)$. We say that the infinitesimal deformation $(\rho_1,\sigma_1)$ is \emph{obstructed to second order} if there exist no second order deformation in its direction.
\end{definition}

Second order deformations coincide with the second derivative of a curve of nearly K\"ahler structures.

\begin{prop}\label{prop:2nd:order:NK:deformation}
Given a nearly K\"ahler structure $(\Real\Omega_0,\tfrac{1}{2}\omega_0^2)$ and an infinitesimal nearly K\"ahler deformation $(\rho_1,\sigma_1)$ suppose that there exists a curve $(\Real\Omega_\epsilon,\tfrac{1}{2}\omega_\epsilon^2)$ of nearly K\"ahler structures such that
\[
(\Real\Omega_\epsilon,\tfrac{1}{2}\omega_\epsilon^2)|_{\epsilon=0}=(\Real\Omega_0,\tfrac{1}{2}\omega_0^2), \qquad \frac{d}{d\epsilon}(\Real\Omega_\epsilon,\tfrac{1}{2}\omega_\epsilon^2)|_{\epsilon=0}=(\rho_1,\sigma_1).
\]
Then $\frac{d^2}{d\epsilon^2}(\Real\Omega_\epsilon,\tfrac{1}{2}\omega_\epsilon^2)|_{\epsilon=0}$ defines a second order deformation in the direction of $(\rho_1,\sigma_1)$.
\end{prop}
In particular, an infinitesimal nearly K\"ahler deformation obstructed to second order cannot be integrated to a curve of nearly K\"ahler structures.

The previous section suggests the following approach to find second (and higher) order deformations of a nearly K\"ahler structure $(\Real\Omega_0,\tfrac{1}{2}\omega_0^2)$. Namely, we look for (formal) power series defining stable exact forms
\[
\begin{gathered}
\Real\Omega_\epsilon = \Real\Omega_0 + \epsilon \rho_1 + \tfrac{1}{2}\rho_2 +\dots,\\
\tfrac{1}{2}\omega_\epsilon ^2 = \tfrac{1}{2}\omega_0^2 + \epsilon \sigma_1 + \tfrac{1}{2}\epsilon^2 \sigma_2 +\dots
\end{gathered}
\]
with $\rho_i\in\Omega^3_{exact}$ and $\sigma_i\in\Omega^4_{exact}$, and a vector field
\[
Z_\epsilon = \epsilon Z_1 + \tfrac{1}{2} \epsilon ^2 Z_2 +\dots, 
\]
satisfying the equation in Proposition \ref{prop:NK:revisited}, \ie
\begin{equation}\label{eq:NK:Iteration}
d\omega_\epsilon -3\Real\Omega_\epsilon =0, \qquad d\Imag\Omega_\epsilon +4\left( \tfrac{1}{2}\omega_\epsilon^2 \right) +d\ast d(Z_\epsilon\Real\Omega_\epsilon)=0,
\end{equation}
where $\Imag\Omega_\epsilon$ and $\omega_\epsilon$ are the duals of $\Real\Omega_\epsilon$ and $\tfrac{1}{2}\omega_\epsilon^2$ respectively.

Given an infinitesimal nearly K\"ahler deformation $(\rho_1,\sigma_1)$, set $Z_1=0$ and look for $(\rho_2,\sigma_2)\in\Omega^3_{exact}\times\Omega^4_{exact}$ such that \eqref{eq:NK:Iteration} are satisfied up to terms of order $O(\epsilon^3)$. Explicitly, we can write the duals $\Imag\Omega_\epsilon$ and $\omega_\epsilon$ of $\Real\Omega_\epsilon$ and $\tfrac{1}{2}\omega_\epsilon^2$, respectively, as
\[
\Imag\Omega_\epsilon = \Imag\Omega_0 + \epsilon \hat{\rho}_1 +\tfrac{1}{2}\epsilon^2 \left( \hat{\rho}_2 - Q_3(\rho_1)\right) +O(\epsilon^3), \qquad \omega_\epsilon = \omega_0 + \epsilon \hat{\sigma}_1 + \tfrac{1}{2}\epsilon^2 \left( \hat{\sigma}_2 - Q_4 (\sigma_1)\right) + O(\epsilon^3),
\]
where $\hat{\, }$ denotes the linearisation of Hitchin's duality map for stable forms in Proposition \ref{prop:Linearisation:Hitchin:dual} and $Q_3(\rho_1)$, $Q_4(\sigma_1)$ are quadratic expressions yielding the quadratic term of Hitchin's duality map. We therefore look for a solution $(\rho_2,\sigma_2,Z_2)$ of
\begin{equation}\label{eq:2nd:order:NK:deformation}
d\hat{\sigma}_2-3\rho_2=dQ_4(\sigma_1), \qquad d\hat{\rho}_2 + 4\sigma_2 +d\ast d(Z_2\lrcorner\Real\Omega_0)=dQ_3 (\rho_1).
\end{equation}
By the previous section we know that there are obstructions to solve these equations. Assume nonetheless that a solution exists. Then we show that $(\rho_2,\sigma_2)$ defines a second order deformation in the sense of Definition \ref{def:2nd:order:NK:deformation}.

\begin{lemma}\label{lem:2nd:order:NK:deformation}
Assume a solution $(\rho_2,\sigma_2,Z_2)$ to \eqref{eq:2nd:order:NK:deformation} exists. Then $d(Z_2\lrcorner\Real\Omega_0)=0$ and
\[
\Big( \rho_2,\hat{\rho}_2-Q_3(\rho_1), \sigma_2, \hat{\sigma}_2-Q_4(\sigma_1)\Big)
\]
defines a second order deformation in the direction of $(\rho_1,\sigma_1)$.
\proof
We need to show that $d(Z_2\lrcorner\Real\Omega_0)=0$ and that the compatibility constraints \eqref{eq:SU(3):structure:Constraints} are satisfied up to order $O(\epsilon^3)$. The proof is similar to the one of Proposition \ref{prop:NK:revisited}.

First observe that since $\omega_\epsilon$ is the dual of $\tfrac{1}{2}\omega_\epsilon^2$, as the notation indicates $\tfrac{1}{2}\omega_\epsilon^2$ is proportional to the square of $\omega_\epsilon$. In particular,
\[
\hat{\sigma}_1\wedge\hat{\sigma}_1 + \omega_0 \wedge \hat{\sigma}_2-\omega_0 \wedge Q_4(\sigma_1)=\sigma_2\in\Omega^4_{exact}.
\]
Hence, using $d\omega_0=3\Real\Omega_0$, $d\hat{\sigma}_1=3\rho_1$ and the first equation in \eqref{eq:2nd:order:NK:deformation} we obtain
\[
3\Big( \hat{\sigma}_1\wedge\rho_1 + \omega_0 \wedge\rho_2 + \hat{\sigma}_2\wedge\Real\Omega_0 - Q_4(\sigma_1)\wedge\Real\Omega_0\Big) = d\Big( \hat{\sigma}_1\wedge\hat{\sigma}_1 + \omega_0 \wedge \hat{\sigma}_2-\omega_0 \wedge Q_4(\sigma_1)\Big)=d\sigma_2=0,
\]
\ie $\omega_\epsilon\wedge\Real\Omega_\epsilon = O(\epsilon^3)$.

In particular, for every vector field $Y$ we have $(Y\lrcorner\Real\Omega_\epsilon)\wedge\tfrac{1}{2}\omega_\epsilon^2 = O(\epsilon^3)$, \ie
\[
8\sigma_1 \wedge (Z\lrcorner\rho_1) + 4\sigma_2 \wedge (Y\lrcorner\Real\Omega_0) + 2\omega_0^2 \wedge (Y\lrcorner\rho_2)=0.
\]

Moreover, by Proposition \ref{prop:Torsion:SU(3):structures} the fact that $d\Real\Omega_\epsilon = O(\epsilon^3)$ implies that $\int_{M}{d\Imag\Omega_\epsilon \wedge (Y\lrcorner \Real\Omega_\epsilon)}=O(\epsilon^3)$ for every vector field $Y$, \ie
\[
\int_{M}{d\Imag\Omega_0 \wedge (Y\lrcorner\rho_2)+ 2d\hat{\rho}_1\wedge (Y\lrcorner\rho_1) + d\left( \hat{\rho}_2-Q_3(\rho_1)\right) \wedge (Y\lrcorner\Real\Omega_0)}=0.
\]
Using $d\Imag\Omega_0=-2\omega_0^2$ and $d\hat{\rho}_1 = -4\sigma_1$ we therefore have
\[
\int_{M}{d\left( \hat{\rho}_2-Q_3(\rho_1)\right) \wedge (Y\lrcorner\Real\Omega_0)}=\int_{M}{8\sigma_1\wedge(Y\lrcorner\rho_1) + 2\omega_0^2\wedge(Y\lrcorner\rho_2)}.
\]

We can now show that $d(Z_2\lrcorner\Real\Omega_0)=0$. Indeed,
\[
\begin{gathered}
\| d(Z_2\lrcorner\Real\Omega_0)\| ^2_{L^2} = -\int_{M}{ (Z\lrcorner\Real\Omega_0)\wedge\Big( d\left( \hat{\rho}_2-Q_3(\rho_1)\right) + 4\sigma_2 \Big)  }\\
=-\int_{M}{ 8\sigma_1 \wedge (Z\lrcorner\rho_1) + 4\sigma_2 \wedge (Y\lrcorner\Real\Omega_0) + 2\omega_0^2 \wedge (Y\lrcorner\rho_2) }=0
\end{gathered}
\]

Finally, since $\Imag\Omega_\epsilon$ is the dual of $\Real\Omega_\epsilon$ and $\omega_\epsilon\wedge\Real\Omega_\epsilon=O(\epsilon^3)$ then $\omega_\epsilon\wedge\Imag\Omega_\epsilon=O(\epsilon^3)$ also, \ie
\[
2\hat{\sigma}_1\wedge\hat{\rho}_1 + \omega_0\wedge \Big( \hat{\rho}_2 -Q_3(\rho_1) \Big) + \Big( \hat{\sigma}_2-Q_4(\sigma_1)\Big)\wedge\Imag\Omega_0=0.
\]
Taking the differential of this expression we find
\[
6\rho_1\wedge\hat{\rho}_1 -8\hat{\sigma}_1\wedge\sigma_1 +3\Real\Omega_0\wedge \Big( \hat{\rho}_2-Q_3(\rho_1)\Big) -4\omega_0\wedge\sigma_2 + 3\rho_2 \wedge \Imag\Omega_0 - 2\Big( \hat{\sigma}_2-Q_4(\sigma_1)\Big)\wedge \omega_0^2=0,
\]
where we used the fact that $d\hat{\rho}_2+4\sigma_2 = dQ_3(\rho_1)$. Up to a constant factor the RHS is exactly the order $\epsilon^2$ coefficient of $3\Real\Omega_\epsilon\wedge\Imag\Omega_\epsilon - 2\omega_\epsilon^3$ and therefore the proof is complete.
\endproof
\end{lemma}

Conversely, every second order deformation $(\rho_2,\rho'_2,\sigma_2,\sigma'_2)$ in the sense of Definition \ref{def:2nd:order:NK:deformation} satisfies $d\sigma'_2-3\rho_2=0$ and $d\rho'_2+4\sigma_2=0$ and the constraints \eqref{eq:SU(3):structure:Constraints} force $\rho'_2 = \hat{\rho}_2-Q_3(\rho_1)$ and $\sigma'_2 = \hat{\sigma}_2-Q_4(\sigma_1)$. Thus we have a complete one-to-one correspondence between second order deformations in the sense of Definition \ref{def:2nd:order:NK:deformation} and solutions to \eqref{eq:2nd:order:NK:deformation}.

We therefore reduced to study the solvability of the equation
\[
d\hat{\sigma}-3\rho=dQ_4(\sigma_1), \qquad d\hat{\rho} + 4\sigma +d\ast d(Z\lrcorner\Real\Omega_0)=dQ_3(\rho_1).
\]
By Proposition \ref{eq:Obstructions} we know that there are obstructions to solve these equations: a solution exists if and only if
\begin{equation}\label{eq:Obstructions}
\langle d^\ast d Q_4(\sigma_1) + 3\ast dQ_3(\rho_1), \hat{\sigma}\rangle _{L^2} = 12\langle Q_4(\sigma_1), \hat{\sigma}\rangle_{L^2} -3\langle Q_3(\rho_1), \ast d\hat{\sigma}\rangle_{L^2}=0
\end{equation}
for every coclosed primitive $(1,1)$ form $\hat{\sigma}$ such that $\triangle\hat{\sigma}=12\hat{\sigma}$.

\subsection{Infinitesimal deformations of the flag manifold}

We now recall the result of Moroianu--Semmelmann \cite[Corollary 6.1]{Moroianu:Semmelmann} on infinitesimal nearly K\"ahler deformations of the flag manifold $\F_3$. We will then show that these infinitesimal deformations are all obstructed to second order.

Introduce left-invariant $1$--forms $e_1,\dots, e_6$ on $\F_3$ dual to the vector fields
\[
E_{12}-E_{21}, \qquad i(E_{12}+E_{21}), \qquad E_{13}-E_{31}, \qquad i(E_{13}+E_{31}), \qquad E_{23}-E_{32}, \qquad i(E_{23}+E_{32}),
\]
where $E_{ij}$ is the $3\times 3$ matrix with $1$ in position $ij$ and all other entries zero.

The complex $1$--forms
\[
\theta_1 = e_2 - ie_1,\qquad \theta_2 = e_4 + i e_3, \qquad \theta_3 = e_6 -ie_5
\]
span the space of $(1,0)$ forms on $\F _3$ with respect to the almost complex structure $J_0$ induced by the homogeneous nearly K\"ahler structure. The nearly K\"ahler structure $(\omega_0,\Omega_0)$ is given by
\[
\omega_0 = \tfrac{i}{2}(\theta_1 \wedge \overline{\theta}_1 + \theta_2\wedge\overline{\theta}_2 + \theta_3 \wedge\overline{\theta}_3), \qquad \Omega_0 = \theta_1\wedge\theta_2\wedge\theta_3.
\]

Given an element $\xi\in\Lie{su}_3$ define functions $v_1,v_2,v_3,x_1,\dots,x_6$ on $\F_3$ by
\[
g^{-1}\xi g = \left( \begin{array}{ccc} \tfrac{i}{2}v_1 & x_1 + i x_2 & x_3 + i x_4\\ -x_1+ix_2 & \tfrac{i}{2}v_2 & x_5 + i x_6\\
-x_3+ix_4 & -x_5 + ix_6 & \tfrac{i}{2}v_3\end{array}\right) = i\left( \begin{array}{ccc} \tfrac{1}{2}v_1 & \overline{z}_1 & z_2\\
z_1 & \tfrac{1}{2}v_2 & \overline{z}_3 \\ \overline{z}_2 & z_3 & \tfrac{1}{2}v_3\end{array}\right),
\]
where
\[
z_1 = x_2 + i x_1, \qquad z_2 = x_4-ix_3, \qquad z_3 = x_6 + ix_5.
\]
Here $g\in SU(3)$ but the functions $v_1,\dots,x_6$ descends to functions on the flag manifold $\F_3=SU(3)/T^2$.

Moroianu--Semmelmann show that for each element $\xi\in\Lie{su}_3$ the $2$--form
\[
\hat{\sigma}_\xi = \tfrac{i}{2}( v_3\, \theta_1 \wedge \overline{\theta}_1 + v_2\, \theta_2\wedge\overline{\theta}_2 + v_1\, \theta_3 \wedge\overline{\theta}_3)
\]
is a primitive coclosed $(1,1)$ form with $\triangle\hat{\sigma}_\xi=12\hat{\sigma}_\xi$.

Define $\rho_\xi$ by $d\hat{\sigma}_\xi=3\rho_\xi$. Since $v_1+v_2+v_3=0$ and $d(\theta_1\wedge\overline{\theta}_1)=d(\theta_2\wedge\overline{\theta}_2)=d(\theta_3\wedge\overline{\theta}_3)$ by \cite[Equation (38)]{Moroianu:Semmelmann}, we have
\[
d\hat{\sigma}_\xi = \tfrac{i}{2}\left( dv_3\wedge \theta_1 \wedge \overline{\theta}_1 + dv_2\wedge \theta_2\wedge\overline{\theta}_2 + dv_1\wedge \theta_3 \wedge\overline{\theta}_3\right).
\]
Moreover, using \cite[Equation (40)]{Moroianu:Semmelmann} and the complex notation we have introduced, we compute
\[
dv_1 = \Imag(z_1\theta_1-z_2\theta_2), \qquad dv_2 = \Imag(z_3\theta_3-z_1\theta_1), \qquad dv_3 = \Imag(z_2\theta_2-z_3\theta_3). 
\] 
Hence
\begin{equation}\label{eq:rho:xi}
\rho_\xi = -\tfrac{1}{6}\Real\left( \theta_1\wedge\theta_2\wedge\alpha_3 +  \theta_2\wedge\theta_3\wedge\alpha_1 + \theta_3\wedge\theta_1\wedge\alpha_2 \right), \qquad \alpha_i = z_j \overline{\theta}_k + z_k \overline{\theta}_j,
\end{equation}
where for each $i=1,2,3$ the two indices $j,k$ are chosen so that $\epsilon_{ijk}=1$.

Set $\sigma_\xi=\hat{\sigma}_\xi\wedge\omega_0$. By Theorem \ref{thm:Infinitesimal:deformations} the pair $(\rho_\xi,\sigma_\xi)$ is an infinitesimal nearly K\"ahler deformation.

\subsection{The quadratic term}

In view of \eqref{eq:Obstructions}, we have now to calculate the quadratic terms $Q_4(\sigma_\xi)$ and $Q_3(\rho_\xi)$.

The $2$--form $Q_4(\sigma_\xi)$ is implicitly defined as follows. Given $\xi$ and chosen $\epsilon>0$ sufficiently small, consider the stable $4$--form $\tfrac{1}{2}\omega_\epsilon = \tfrac{1}{2}\omega_0^2 + \epsilon\sigma_\xi$ and its dual
\[
\omega_\epsilon=\omega_0 + \epsilon \hat{\sigma}_\xi - \tfrac{\epsilon^2}{2}Q_4(\hat{\sigma}_\xi) + O(\epsilon^3).
\]
As the notation suggests $\tfrac{1}{2}\omega_\epsilon^2$ is proportional to the square of $\omega_\epsilon$. In particular,
\begin{equation}\label{eq:Quadratic:4:form}
\hat{\sigma}_\xi\wedge\hat{\sigma}_\xi-\omega_0 \wedge Q_4 (\sigma_\xi)=0.
\end{equation}
From this formula we could calculate explicitly $Q_4(\sigma_\xi)$ but we will see below that this is not necessary.

Similarly, $Q_3 (\rho_\xi)$ is defined so that
\[
\Imag\Omega_0 - \epsilon\ast\rho_\xi -\tfrac{\epsilon^2}{2}Q_3(\rho_\xi) +O(\epsilon^3) 
\]
is the dual of the stable $3$--form $\Real\Omega_0 + \epsilon \rho_\xi$. While it is possible to compute $Q_3(\rho_\xi)$ explicitly following the algorithm to construct the dual of a stable $3$--form \cite[Equation (2) and \S 8.2]{Hitchin} we find it quicker to produce by hand a complex volume form $\Omega_\epsilon$ such that
\[
\Real\Omega_\epsilon = \Real\Omega_0 + \epsilon \rho_\xi + O(\epsilon^3).
\]
Then the quadratic term of $\Imag\Omega_\epsilon$ will yield an explicit expression for $Q_3(\rho_\xi)$.

Consider the complex $1$--forms
\[
\theta_i^\epsilon = \theta_i -\tfrac{\epsilon}{6}(z_k\overline{\theta}_j + z_j \overline{\theta}_k) + \tfrac{\epsilon^2}{36}\,\overline{z}_i\, (\overline{z}_i \overline{\theta}_i -\overline{z}_j \overline{\theta}_j -\overline{z}_k \overline{\theta}_k ),
\]
where as before for each $i=1,2,3$ the two indices $j,k$ are chosen so that $\epsilon_{ijk}=1$. For $\epsilon$ sufficiently small $\theta_1^\epsilon, \theta_2^\epsilon, \theta_3^\epsilon$ span a $3$--dimensional subspace of the space of complex $1$--forms and therefore define an almost complex structure. We now set $\Omega_\epsilon = \theta_1^\epsilon\wedge \theta_2^\epsilon\wedge \theta_3^\epsilon$.

Expanding $\Omega_\epsilon = \Omega_0 + \epsilon \Omega_1 + \epsilon^2 \Omega_2 + O(\epsilon^3)$, one can check that $\rho_\xi = \Real\Omega_1$ and $\Real\Omega_2=0$. Thus $Q_3(\rho_\xi)=-2\Imag\Omega_2$. Via a straightforward calculation we conclude that
\begin{equation}\label{eq:Quadratic:3:form}
9\, Q_3(\rho_\xi) = \Imag\left( \overline{\theta}_1\wedge\overline{\theta}_2\wedge\beta_3 + \overline{\theta}_2\wedge\overline{\theta}_3\wedge\beta_1 + \overline{\theta}_3\wedge\overline{\theta}_1\wedge\beta_2 \right), \qquad \beta_i = z_i(z_i \theta_i - z_j\theta_j - z_k\theta_k), 
\end{equation}
with the usual convention for indices.

\subsection{Representation theoretic considerations}

We now use basic representation theoretic observations to reduce to a minimum the number of calculations necessary to prove that infinitesimal deformations of the flag manifold are obstructed to second order.

By \eqref{eq:Obstructions} for every $\xi\in\Lie{su}_3$ the infinitesimal nearly K\"ahler deformation $(\rho_\xi,\sigma_\xi)$ is integrable to second order if and only if
\[
\Phi(\xi,\xi'):= \langle 12\, Q_4(\sigma_\xi),\hat{\sigma}_{\xi'}\rangle_{L^2} -3\langle Q_3(\rho_\xi),\ast d\hat{\sigma}_{\xi'}\rangle_{L^2}=0
\]
for every $\xi'\in\Lie{su}_3$.

By invariance of the nearly K\"ahler structure with respect to the natural action of $SU(3)$ on $\F_3$, $\Phi$ defines an $\text{Ad}$--invariant map $\Phi\co \Lie{su}_3\times\Lie{su}_3\ra \R$ which is quadratic in the first argument and linear in the second. It must therefore correspond to an element of $\text{Hom}_{SU(3)}(\Lie{su}_3,Sym^2(\Lie{su}_3))$. This is a one dimensional space. A generator $\xi\mapsto L_\xi$ can be defined as follows. Identify $\Lie{su}_3$ with $\Lambda^2_8(\R^6)^\ast$. Given $\xi,\xi'\in\Lie{su}_3$ we have
\[
\xi\wedge\xi'=-\tfrac{1}{6}\langle \xi,\xi'\rangle \omega_0 + \omega_0 \wedge L_{\xi}(\xi')
\]
for a unique primitive $(1,1)$ form $L_{\xi}(\xi')$. Here $\omega_0$ is the standard $SU(3)$--invariant K\"ahler form on $\R^6\simeq \C^3$. We therefore conclude that $\Phi (\xi,\xi')=C \Phi_0 (\xi^2,\xi')$ for some constant $C\in \R$. Here $\xi^2\in Sym^2(\Lie{su}_3)$ and 
\[
\Phi_0 (\xi^2,\xi') = \langle L_{\xi'}, \langle \xi, \cdot \rangle\xi \rangle.
\] 

Furthermore, we will now argue that in order to determine whether $C\neq 0$ it is enough to calculate the numbers $\Phi (\xi,\xi)$ for $\xi\in\Lie{su}_3$. To this end, it is enough to show that $\Phi_0(\xi^2,\xi)$ is a non-zero multiple of the unique cubic invariant polynomial on $\Lie{su}_3$:
\[
\Phi_0(\xi^2,\xi)=\langle L_\xi, \langle \xi, \cdot \rangle\xi \rangle = \sum_{i=1}^8{\langle L_\xi(e_i), \langle \xi, e_i \rangle\xi \rangle} = \langle L_\xi (\xi) ,\xi \rangle = \xi\wedge\ast L_\xi (\xi)=-\xi\wedge\xi\wedge\xi
\]
by Lemma \ref{lem:Hodge:star}.(iii). Here $e_1,\dots, e_8$ is an orthonormal basis of $\Lie{su}_3$.

\subsection{Obstructed deformations}

We have now all the ingredients to prove the main theorem of this section.

\begin{theorem}\label{thm:Flag}
The infinitesimal deformations of the homogeneous nearly K\"ahler structure on the flag manifold $\F_3$ are all obstructed.
\proof
By Proposition \ref{prop:2nd:order:NK:deformation} it is enough to  show that the infinitesimal deformations are obstructed to second order in the sense of Definition \ref{def:2nd:order:NK:deformation}.

Consider the map $\Phi\co \Lie{su}_3\times\Lie{su}_3\ra \R$ introduced in the previous section and in particular its values on vectors of the form $(\xi,\xi)$ with $\xi\in\Lie{su}_3$. We know that $\Phi(\xi,\xi)=iC \det(\xi)$ for some $C\in\R$. We want to show that $C\neq 0$.

Now, $\Phi(\xi,\xi)=\int_{\F_3}{f_\xi ([g]) \dvol}$, where $[g]\in \F_3 = SU(3)/T^2$ and $f_\xi$ is the function on $\F_3$ defined by
\[
f_\xi ([g]) = \langle 12 Q_4(\sigma_\xi),\hat{\sigma}_{\xi}\rangle -3\langle Q_3(\rho_\xi),\ast d\hat{\sigma}_{\xi}\rangle.
\]

By the definition of $\hat{\sigma}_\xi$ and $3\rho_\xi=d\hat{\sigma}_\xi$ in terms of the functions $v_1,\dots, x_6$, it is clear that $f_\xi ([h^{-1}g])=f_{h\xi h^{-1}}([g])$. We therefore conclude that $f_\xi([g])=P(g^{-1}\xi g)$, where $P$ is a polynomial of degree $3$ on $\Lie{su}_3$.

In fact we can use our formulas for $Q_4(\sigma_\xi)$ and $Q_3(\rho_\xi)$ to explicitly calculate $P$. Indeed, using $d\hat{\sigma}_\xi=3\rho_\xi$ we have 
\[
f_\xi ([g])\dvol =12\, Q_4(\sigma_\xi)\wedge \ast\hat{\sigma}_\xi +9\, Q_3(\rho_\xi)\wedge\rho_\xi.  
\]
Now use \eqref{eq:Quadratic:4:form} and the fact that $\hat{\sigma}_\xi\in\Omega^2_8$ to write
\[
12\, Q_4(\sigma_\xi)\wedge \ast\hat{\sigma}_\xi = -12\, Q_4(\sigma_\xi)\wedge\omega_0\wedge \hat{\sigma}_\xi = -12\,\hat{\sigma}_\xi\wedge \hat{\sigma}_\xi\wedge \hat{\sigma}_\xi = -12\, v_1 v_2 v_3 \, \omega_0^3.
\]
On the other hand, using \eqref{eq:rho:xi} and \eqref{eq:Quadratic:3:form} we calculate
\[
\begin{gathered}
9\, Q_3(\rho_\xi)\wedge\rho_\xi = -\tfrac{1}{6}\Imag\left( \sum_{i=1}^3{\overline{\theta}_j\wedge\overline{\theta}_k\wedge\beta_i}\right)\wedge \Real\left( \sum_{i=1}^3{\theta_j\wedge\theta_k\wedge\alpha_i}\right)=\\
-\tfrac{1}{12}\Imag\left( ( \overline{\theta}_1\wedge\overline{\theta}_2\wedge\beta_3 + \overline{\theta}_2\wedge\overline{\theta}_3\wedge\beta_1 + \overline{\theta}_3\wedge\overline{\theta}_1\wedge\beta_2 )\wedge \left( \theta_1\wedge\theta_2\wedge\alpha_3 +  \theta_2\wedge\theta_3\wedge\alpha_1 + \theta_3\wedge\theta_1\wedge\alpha_2 \right)\right)\\
= \tfrac{1}{2}\Imag\left( z_1 z_2 z_3\, \theta_1\wedge\theta_2\wedge\theta_3\wedge \overline{\theta}_1\wedge\overline{\theta}_2\wedge\overline{\theta}_3\right)=-\Real(z_1 z_2 z_3) \Real\Omega_0 \wedge\Imag\Omega_0,
\end{gathered}
\]
since $\Omega_0=\theta_1\wedge\theta_2\wedge\theta_3$ and $\Omega_0\wedge\overline{\Omega}_0=-2i\Real\Omega_0\wedge\Imag\Omega_0$.

We conclude that
\[
f_\xi([g]) =P(g^{-1}\xi g)= -72\, v_1 v_2 v_3 - 4\Real\left( z_1 z_2 z_3\right).
\]

In order to calculate the mean value of $P$ on $\mathbb{F}_3$ we appeal to the Peter--Weyl Theorem.

First of all observe that we can lift $f_\xi$ to $SU(3)$ as a $T^2$--invariant function and calculate $\Phi(\xi,\xi)$ up to a positive factor by considering the mean value of $f_\xi$ on $SU(3)$. Indeed, by \cite[Lemma 5.4]{Moroianu:Semmelmann} the nearly K\"ahler metric on $\F_3$ is induced by $-\tfrac{1}{12}B$, where $B$ is the Killing form of $SU(3)$.

The Peter-Weyl Theorem says that for any compact Lie group $G$
\[
L^{2}(G) = \bigoplus _{\gamma \in \hat{G}}{ V_{\gamma} \otimes V_{\gamma}^{\ast}},
\]
where $\hat{G}$ is the set of (non-isomorphic) irreducible $G$--representations. The Peter--Weyl isomorphism is explicit: to a pair $v \otimes \alpha \in V _{\gamma} \otimes V_{\gamma}^\ast$ we associate the function $f(g)=\alpha \left( \gamma (g^{-1})v \right)$. Moreover, each summand of the Peter--Weyl decomposition is an eigenspace for the Laplacian $\triangle _{G}$ with eigenvalue that can be determined from the highest weight $\gamma$.

Now, the decomposition of $Sym^3(\Lie{su}_3)$ into irreducible representations of $SU(3)$ contains a unique copy of the trivial representation, corresponding to the unique cubic $\text{Ad}$--invariant polynomial $i\det$ on $\Lie{su}_3$. Thus we write $Sym^3(\Lie{su}_3)=V_0\oplus\dots \oplus V_n$ with $V_0\simeq \R$ and $V_i$, $i>0$, non trivial representations $\gamma_i$ (it is not important what these actually are nor that each of these appears without multiplicity). We must have $P(g^{-1}\xi g)=\sum_{i=0}^n{\alpha_i \left( \gamma_i(g^{-1})(\xi^3)_i\right)}$ for some $\alpha_i\in V_i^\ast$. Here for every $\xi\in\Lie{su}_3$ we write $\xi^3=(\xi^3)_0 +\dots + (\xi^3)_n$ according to the decomposition of $Sym^3(\Lie{su}_3)$. In this way we obtain the decomposition of $f_\xi$ into eigenspaces for the Laplacian on $SU(3)$. In particular, we can compute the mean value of $f_\xi$ simply by calculating the inner product of the cubic polynomial $P$ with the unique invariant cubic polynomial $i\det$.

Now, $i\det (g^{-1}\xi g)=\tfrac{1}{8}v_1 v_2 v_3 -\tfrac{1}{2}(v_3 |z_1|^2 + v_2 |z_2|^2 + v_1 |z_3|^2) +2\Real(z_1 z_2 z_3)$. Since the monomials $v_1 v_2 v_3$ and $\Real (z_1 z_2 z_3)$ appear with coefficients of the same sign both in $P$ and $i\det$ and $\text{Span}(e_1,\dots,e_6)$ is orthogonal to the sub-algebra of diagonal matrices (so that $v_1 v_2 v_3$ and $\Real (z_1 z_2 z_3)$ are orthogonal polynomials), we conclude that $P$ has a non-zero inner product with $i\det$ and the proof is complete.
\endproof 
\end{theorem}

\bibliographystyle{amsinitial}
\bibliography{nkdef}

\providecommand{\bysame}{\leavevmode\hbox to3em{\hrulefill}\thinspace}
\providecommand{\MR}{\relax\ifhmode\unskip\space\fi MR }
\providecommand{\MRhref}[2]{%
  \href{http://www.ams.org/mathscinet-getitem?mr=#1}{#2}
}
\providecommand{\href}[2]{#2}
\begin{thebibliography}{10}

\bibitem{Alexandrov:Semmelmann}
B.~Alexandrov and U.~Semmelmann, \emph{Deformations of nearly parallel {${\rm
  G}_2$}-structures}, Asian J. Math. \textbf{16} (2012), no.~4, 713--744.

\bibitem{Bar}
C.~B{\"a}r, \emph{Real {K}illing spinors and holonomy}, Comm. Math. Phys.
  \textbf{154} (1993), no.~3, 509--521.

\bibitem{Friedrich:al}
H.~Baum, T.~Friedrich, R.~Grunewald, and I.~Kath, \emph{Twistors and {K}illing
  spinors on {R}iemannian manifolds}, Teubner-Texte zur Mathematik [Teubner
  Texts in Mathematics], vol. 124, B. G. Teubner Verlagsgesellschaft mbH,
  Stuttgart, 1991, With German, French and Russian summaries.

\bibitem{Bedulli:Vezzoni:SU(3)}
L.~Bedulli and L.~Vezzoni, \emph{The {R}icci tensor of {SU}(3)-manifolds}, J.
  Geom. Phys. \textbf{57} (2007), no.~4, 1125--1146.

\bibitem{Charbonneau:Harland}
B.~Charbonneau and D.~Harland, \emph{{D}eformations of nearly {K}{\"a}hler
  instantons}, 2015, arXiv:1510.07720.

\bibitem{Cortes}
V.~Cort{\'e}s and J.~J. V{\'a}squez, \emph{Locally homogeneous nearly
  {K}\"ahler manifolds}, Ann. Global Anal. Geom. \textbf{48} (2015), no.~3,
  269--294.

\bibitem{Foscolo:Haskins}
L.~Foscolo and M.~Haskins, \emph{{N}ew {G}2 holonomy cones and exotic nearly
  {K}aehler structures on the 6-sphere and the product of two 3-spheres}, 2015,
  arXiv:1501.07838.

\bibitem{Friedrich}
T.~Friedrich, \emph{Der erste {E}igenwert des {D}irac-{O}perators einer
  kompakten, {R}iemannschen {M}annigfaltigkeit nichtnegativer
  {S}kalarkr\"ummung}, Math. Nachr. \textbf{97} (1980), 117--146.

\bibitem{Goto}
R.~Goto, \emph{Moduli spaces of topological calibrations, {C}alabi-{Y}au,
  hyper-{K}\"ahler, {$G_2$} and {${\rm Spin}(7)$} structures}, Internat. J.
  Math. \textbf{15} (2004), no.~3, 211--257.

\bibitem{Gray}
A.~Gray, \emph{Nearly {K}\"ahler manifolds}, J. Differential Geometry
  \textbf{4} (1970), 283--309.

\bibitem{Hitchin}
N.~Hitchin, \emph{Stable forms and special metrics}, Global differential
  geometry: the mathematical legacy of {A}lfred {G}ray ({B}ilbao, 2000),
  Contemp. Math., vol. 288, Amer. Math. Soc., Providence, RI, 2001, pp.~70--89.

\bibitem{Karigiannis:Lotay}
S.~Karigiannis and J.~Lotay, \emph{{D}eformation theory of {G}2 conifolds},
  2012, arXiv:1212.6457.

\bibitem{Koiso}
N.~Koiso, \emph{Rigidity and infinitesimal deformability of {E}instein
  metrics}, Osaka J. Math. \textbf{19} (1982), no.~3, 643--668.

\bibitem{Moroianu:Semmelmann:Unit:Killing}
A.~Moroianu, P.-A. Nagy, and U.~Semmelmann, \emph{Unit {K}illing vector fields
  on nearly {K}\"ahler manifolds}, Internat. J. Math. \textbf{16} (2005),
  no.~3, 281--301.

\bibitem{Moroianu:Nagy:Semmelmann}
\bysame, \emph{Deformations of nearly {K}\"ahler structures}, Pacific J. Math.
  \textbf{235} (2008), no.~1, 57--72.

\bibitem{Moroianu:Semmelmann}
A.~Moroianu and U.~Semmelmann, \emph{The {H}ermitian {L}aplace operator on
  nearly {K}\"ahler manifolds}, Comm. Math. Phys. \textbf{294} (2010), no.~1,
  251--272.

\bibitem{Nordstrom:Thesis}
J.~Nordstr{\"o}m, \emph{Deformations and gluing of asymptotically cylindrical
  manifolds with exceptional holonomy}, Ph.D. thesis, University of Cambridge,
  Cambridge, 2008.

\bibitem{Nordstrom:ACyl}
\bysame, \emph{Deformations of asymptotically cylindrical {$G_2$}-manifolds},
  Math. Proc. Cambridge Philos. Soc. \textbf{145} (2008), no.~2, 311--348.

\bibitem{Pedersen:Poon}
H.~Pedersen and Y.~S. Poon, \emph{A note on rigidity of {$3$}-{S}asakian
  manifolds}, Proc. Amer. Math. Soc. \textbf{127} (1999), no.~10, 3027--3034.

\bibitem{Tipler:vanCoevering}
C.~van Coevering and C.~Tipler, \emph{Deformations of {C}onstant {S}calar
  {C}urvature {S}asakian {M}etrics and {K}-{S}tability}, Int. Math. Res.
  Notices \textbf{2015} (2015), no.~22, 11566--11604.

\bibitem{Verbitsky}
M.~Verbitsky, \emph{Hodge theory on nearly {K}\"ahler manifolds}, Geom. Topol.
  \textbf{15} (2011), no.~4, 2111--2133.

\bibitem{Gray:Wolf}
J.~A. Wolf and A.~Gray, \emph{Homogeneous spaces defined by {L}ie group
  automorphisms. {II}}, J. Differential Geometry \textbf{2} (1968), 115--159.

\end{thebibliography}

\end{document}